\documentclass[11pt]{amsart}

\usepackage{amscd}
\usepackage{amssymb}
\usepackage{boxedminipage}
\usepackage{enumerate}

\theoremstyle{plain} \newtheorem{theorem}{Theorem}[section]
\theoremstyle{plain} \newtheorem{lemma}[theorem]{Lemma}
\theoremstyle{plain} \newtheorem{proposition}[theorem]{Proposition}

\newtheorem{corollary}[theorem]{Corollary}

\newcommand{\nr}{\refstepcounter{theorem}  
                   \noindent {\thetheorem .}}

\newcommand{\eks}{\medskip \noindent {\it Example \nr} }
\newcommand{\eksfin}{\medskip}
\newcommand{\rem}{\medskip \noindent {\it Remark \nr} }
\newcommand{\remfin}{\medskip}

\theoremstyle{definition}

\newcommand{\Hom}{\text{Hom}}

\newcommand{\Ext}{\text{Ext}}

\newcommand{\sus}{\subseteq}
\newcommand{\pil}{\rightarrow}
\newcommand{\vpil}{\leftarrow}

\newcommand{\inpil}{\hookrightarrow}

\newcommand{\mto}[1]{\stackrel{#1}\longrightarrow}

\newcommand{\iso}{\cong}
\newcommand{\te}{\otimes}

\newcommand{\gr}{\text{gr}}

\newcommand{\Ip}{I^\prime}

\newcommand{\ga}{\gamma}

\newcommand{\bev}{B^{ev}}
\newcommand{\bod}{B^{od}}
\newcommand{\one}{{\bf 1}}
\newcommand{\ainf}{A_{\infty}}
\newcommand{\mlto}[1]{\stackrel{#1}\longleftarrow}

\def\hpil{\longrightarrow}

\def\surj{\hpil\mspace{-26.0mu}\hpil}

\newcommand{\ent}[1]{{\bf 1}^{\te{#1}}}

\newcommand{\en}{{\bf 1}}
\newcommand{\smpm}{{\scriptstyle{\pm}}}

\begin{document}

\title [PBW-deformations of N-Koszul algebras]
{PBW-deformations of N-Koszul algebras}
\author { Gunnar Fl{\o}ystad \and Jon Eivind Vatne}
\address{ Matematisk Institutt\\
          Johs. Brunsgt. 12 \\
          5008 Bergen \\
          Norway}   
        
\email{ gunnar@mi.uib.no \and jonev@mi.uib.no}

\begin{abstract} 
For a quotient algebra $U$ of the tensor algebra we give explicit
conditions on its relations for $U$ being a PBW-deformation of an
$N$-Koszul algebra $A$. We show there is a one-to-one correspondence
between such deformations and a class of $A_\infty$-structures
on the Yoneda algebra $Ext_A^*(k,k)$ of $A$. We compute the
PBW-deformations of the algebra whose relations are the anti-symmetrizers
of degree $N$ and also of cubic Artin-Schelter algebras.
\end{abstract}


\maketitle

\section*{Introduction}

In this paper we consider a class of deformations of $N$-Koszul
algebras, the Poincar\'{e}-Birkhoff-Witt(PBW)-deformations.  
The class of 
$N$-Koszul algebras was introduced by R. Berger in \cite{Ber01} and
is a generalization of ordinary Koszul algebras.
Let $R$ be a subspace of homogeneous forms of degree $N$
in the tensor algebra $T(V)$. Denoting by $(R)$ the two-sided 
ideal in the tensor algebra generated by $R$ we get the quotient algebra
$A = T(V)/(R)$.
When certain nice homological conditions are imposed on the resolution
of the residue field $k$ this is called 
an $N$-Koszul algebra, 
the quadratic case $N=2$ giving ordinary Koszul algebras.

By deforming $R$ to a space $P$ of non-homogeneous forms of degree $N$
but keeping the homogeneous part of $P$ of degree $N$ fixed equal to $R$,
we get a filtered algebra $U = T(V)/(P)$ which is a {\it PBW-deformation} of
$A$ if its associated graded algebra is $A$. (Note that $U$ is augmented
if the constant term in $P$ is zero.) In the quadratic case such
deformations have been considered by Braverman and Gaitsgory in \cite{BG96},
a typical example being when $A$ is the symmetric algebra on a vector space $V$
and $U$ the enveloping algebra of a Lie algebra on $V$. We generalize
their results to the setting of $N$-Koszul algebras and give a simple
set of equations relating the homogeneous parts of $P$ which describes when
$U$ gives a PBW-deformation.

In \cite{Po}, Positselski in the quadratic case showed that giving an
augmented PBW-deformation of $A$ is equivalent to giving the Koszul dual 
$A^!$ of $A$ the structure of a differential graded algebra.
The Koszul dual $A^!$ may be identified as the Yoneda
algebra $\Ext^*_A(k,k)$. 
In general for any augmented algebra the Yoneda algebra has the 
structure of an $A_\infty$-algebra. For an $N$-Koszul algebra $A$
it is natural to consider its Yoneda algebra $B$ together with this
$A_\infty$-structure as its Koszul dual and Keller \cite{KeNK}
computed this. We show that when we specify a certain class of 
properties of $A_\infty$-structures of $B$, 
among them being that those 
products which are non-zero in the case computed by Keller are 
kept fixed,
there is a one-to-one correspondence between such $A_\infty$-structures
on $B$ and augmented PBW-deformations of $A$. 
Our description
is quite explicit. In fact the properties we can impose on $B$ reduce
the whole menagerie of equations to be verified for an $A_\infty$-algebra
to only two simple types.
We also verify that the $A_\infty$-structure on the Yoneda algebra
of $U$ corresponds to, more precisely is quasi-isomorphic to, 
the $A_\infty$-structure on $B$ which we find.

\medskip

We then proceed to compute PBW-deformations of some classes of concrete
algebras. The most well known classical case is of course the PBW-deformations
of the symmetric algebra $S(V)$. Such a deformation $U$ corresponds
to giving $V$ the structure of a Lie algebra and  $U$ is then the enveloping
algebra of $V$. Here we consider the class of algebras
\begin{equation}
T(V)/(\wedge^N V)\label{InLikSN}
\end{equation}
whose $N$-Koszulity was shown by Berger \cite{Ber01}. Note that when $N=2$
we get the symmetric algebra $S(V)$. We compute the PBW-deformations
of these algebras. The result separates into two cases according to
$N$ being even or odd. When $N$ is even, an augmented PBW-deformation of 
(\ref{InLikSN})
corresponds to giving a Lie algebra structure $L : \wedge^2 V \pil V$ on
$V$ together with additional structure provided by even symplectic forms
$\Phi_{2r} : \wedge^{2r} V \pil k$ where $2 \leq 2r < N$, fulfilling
analogues of the Jacobi identity, namely the composition (with some abuse
of notation)
\[L \circ \Phi_{2r} \circ L : \wedge^{2r+3} V \pil V \]
is zero. Note that if we take $\Phi_0 = 1$ this is the ordinary
Jacobi identity.

When $N$ is odd a PBW-deformation of (\ref{InLikSN}) is equivalent
to freely choosing a linear form $l : V \pil k$ and symplectic forms
$\Phi_{2r}$ as above. No relations among them is required.

\medskip

We also consider the (generic forms of)  
cubic Artin-Schelter algebras and compute their PBW-deformations.
These are algebras that have considerable interest in non-commutative
algebraic geometry.

\medskip

We would like to give a quick comment about an aspect we have not
discussed in this paper.  Deformations of relations can be
thought of as changes to the Massey product structure.  With this in
mind, the results about the $A_\infty$-structure can be seen as an
explicit example of the interplay between Massey products and
$\infty$-structures.  Also, the
space of deformations we construct will be the versal family of the
moduli of such deformations.  This means that we don't consider
equivalences of deformations, though the quotient by this equivalence
relation gives the actual moduli.  The moduli is in general difficult
to understand, even in seemingly simple situations.  For the case of
Lie algebras on a three-dimensional vector space (corresponding to
PBW-deformations of $k[x,y,z]$), see \cite{Siq04}.

\medskip

The organization of the paper is as follows. Section 1 introduces 
$N$-Koszul algebras and their PBW-deformations. We give the PBW-theorem
for $N$-Koszul algebras which shows which equations  must be fulfilled
for the relations $P$. In Section 2 we state how a PBW-deformation
of $A$ corresponds to an $A_\infty$-structure on its Yoneda algebra
and Section 3 contains the main bulk of the proof of this.

Sections 4, 5, 6, and 7 are concerned with the PBW-deformations of the algebra
(\ref{InLikSN}). In Section 4 we state the results. Sections 5 and 6 are of
a technical nature and develop most of the ingredients of
the proofs. The final proofs are given in Section 7.

Section 8 contains explicit computations of the PBW-deformations of 
cubic Artin-Schelter algebras.

\medskip

\noindent {\bf Acknowledgements}
We would like to thank the Mittag-Leffler institute, which the first
author visited during parts of, and the second
author for the majority of this work.  We would also like to thank
A. Laudal for interesting conversations on this subject.

\section{$N$-Koszul algebras and the PBW-theorem}

\subsection{Definitions}

Let $V$ be a finite dimensional vector space, whose dimension 
we shall denote by $v$, over a field $k$. 
Let $R$ be a subspace of $V^{\te N}$ where $N$ is some
integer $\geq 2$, and $(R)$ the two-sided ideal generated by $R$ 
in the tensor algebra
$T(V)$. We then get a homogeneous algebra $A = T(V)/(R)$.
In \cite{Ber01}, R. Berger introduced the notion of $N$-Koszul algebras
which is a generalization of ordinary Koszul algebras (which is the
case $N=2$). The algebra $A$ is an {\it $N$-Koszul algebra} if  the residue
field $k$ has a right $A$-module resolution of the form 
\[ k \vpil A \vpil V_1 \te A(-1) \vpil V_2 \te A(-N) \vpil \cdots
\vpil V_i \te A(-n(i)) \vpil \cdots \]  
where $n(2q) = Nq$ and $n(2q+1) = Nq+1$. The $V_i$ are vector spaces
and here and in the rest of the paper all tensor products are over
$k$.  We can also define this notion using a resolution of $k$ as a 
left module; these two definitions are then equivalent for a given
algebra, see \cite{Ber01} Definition 2.10, and the following discussion.

 Instead of a homogeneous space of relations $R$ we may consider 
a non-homogeneous space of relations $P$ in $\oplus_{i \leq N} V^{\te i}$,
and get a non-homogeneous algebra $U = T(V)/(P)$. We shall be interested
in the case when this algebra is a deformation of $A$ of a particular kind, a
PBW-deformation. We then assume that $P$ intersects 
$\oplus_{i \leq N-1} V^{\te i}$
trivially and that $R$ is the image of $P$ by the natural projection 
of $\oplus_{i \leq N} V^{\te i}$ to $V^{\te N}$. Then there are linear
maps $R \mto{\alpha_i} V^{\te N-i}$ such that we may write
\begin{equation} 
P = \{ r + \alpha_1(r) + \cdots + \alpha_N(r) \, | \, r \in R \}. 
\label{NKLikP}\end{equation}

There is a natural filtration on the tensor algebra by letting
$F^lT(V)$ be $\oplus_{i \leq l} V^{\te i}$.  This induces a natural 
filtration $F^lU$ on the quotient algebra $U$, and in the
situation described there is a surjection 
\[ A \surj \text{gr}\,U. \]
We say that $U$ is a {\it PBW-deformation} of $A$ if this map is an 
isomorphism. 

\medskip
\noindent {\bf Notation.} The symbol $\en$ will be used to denote
the identity map on a vector space. If $\phi$ and $\psi$ are
maps defined on subspaces $U$ resp. $W$ of $V^{\te a}$ resp.
$V^{\te b}$, denote
\[ [\phi, \psi] = \phi \te \psi - \psi \te \phi \]  
which is defined on the subspace $(U\te W) \cap (W \te U)$ of 
$V^{\te a+b}$.

\subsection{The PBW-theorem}
The topic of this paper is to investigate PBW-deformations of 
$N$-Koszul algebras. In the quadratic case this has been done
in various papers, \cite{BG96}, \cite{Po}, \cite{PP94}.
The following is a generalization of the main theorem of \cite{BG96}
to the $N$-Koszul case.

\begin{theorem} \label{NKThmPBW}
Let $A$ be an $N$-Koszul algebra. Then $U$ defined by relations
(\ref{NKLikP}) is a PBW-deformation of $A$ if and only if 
\begin{equation} 
(P) \cap F^NT(V) = P. \label{NKLikPF}
\end{equation} 
This is more explicitly equivalent to the image of 
\begin{equation} 
(V\te R) \cap (R \te V) \mto{[\en, \alpha_1 ]} V^{\te N} \label{NKLikJac1}
\end{equation}
being contained in $R$ and (let $\alpha_{N+1} = 0$)
\begin{equation}
\alpha_i \circ [\en, \alpha_1 ] = [\en, \alpha_{i+1}], \quad
i = 1, \ldots, N 
 \label{NKLikJac2}
\end{equation}
\end{theorem}

\eks When the relations $R$ are $\wedge^2 V$ so $\alpha_1$ is 
a map $\wedge^2 V \pil V$, and $U$ is augmented, i.e. 
$\alpha_2$ is zero, there is only one equation (\ref{NKLikJac2})
and written out this is just the Jacobi identity meaning that 
$\alpha_1$ gives $V$ the structure of a Lie algebra. Thus we have
the classical PBW-theorem that $U$ is a PBW-deformation of the 
symmetric algebra $S(V)$ iff $U$ is the enveloping algebra
of a Lie algebra.
\eksfin

\rem When $U$ is a PBW-deformation (\ref{NKLikPF}), (\ref{NKLikJac1})
and (\ref{NKLikJac2}) still hold without
the assumption that $A$ is $N$-Koszul, as we now show in the first part
of the proof. 
\remfin

\rem Berger and Ginzburg, \cite{BG05}, have independently proved this result.
\remfin 

\begin{proof}[Proof of Thm. \ref{NKThmPBW}, first part.]
When $U$ is a PBW-deformation of $A$, clearly $(P) \cap F^NT(V)$ is 
equal to $P$.
 
Let $x$ be in 
$(V\otimes R) \cap (R\otimes V)$.
We get in the ideal $(P)$ two elements
\[ \sum_{i=1}^{N} (\en \otimes \alpha_i) (x) + x, \quad
\sum_{i=1}^{N}(\alpha_i\otimes \en) (x) + x \]
The difference of these two elements is in $(P)\cap F^NT(V)$ which is $P$.
This gives us the necessary relations degree by degree.  For
instance, in degree $N$ we have
\[[\en,\alpha_1](x)\]
so this must be in $R$, and the degree $N-1$ part must be $\alpha_1$
of this element, which gives us
\[\alpha_1\circ [\en,\alpha_1](x)= 
[\en,\alpha_2](x)\]
and so forth.

\end{proof}

The second and difficult part of the proof will show how the conditions 
(\ref{NKLikJac1}) and (\ref{NKLikJac2}) of Theorem \ref{NKThmPBW} imply that 
$U$ is a PBW-deformation 
of $A$. For this we need some deformation theory.

\subsection{Deformations}
We recall the setup from \cite{BG96}.
Let $A=\bigoplus_{d=0}^{\infty}A_d$ be a graded
associative algebra over a field $k$.  An $i$th level graded
deformation of $A$ is a graded $k[t]/(t^{i+1})$-algebra $A^i$, free as
a module over $k[t]/(t^{i+1})$, 
together with an isomorphism $A^i/tA^i\iso A$. By taking
inverse limits, a graded deformation of $A$ is a graded
$k[[t]]$-algebra $A^t$, free as a module over $k[[t]]$, with an isomorphism
$A^t/tA^t\iso A$.  In all cases considered in this paper, the
deformations we construct are naturally defined also over $k[t]$, with
the same conditions.  Thus we can, for instance, take fibres over the
point $t=1$.

More explicitly, an $i$th level deformation is given by a set of maps
$f_j:A\otimes A\rightarrow A$ such that the rule
\[a\times b=ab+\sum_{j=1}^{i}f_j(a,b)t^j\]
defines an associative algebra, and similarly for the inverse limit,
but with the sum extended to infinity.\\

Let
$\mathcal{C}, \mathcal{C}^{i}$ denote the categories of graded
deformations and graded $i$th level deformations of $A$, respectively (with
morphisms only isomorphisms).

The connection between graded Hochschild cohomology groups and
deformations is given by the next lemma, taken directly from
\cite{BG96}.

\begin{lemma}\label{HHbet}
\begin{itemize}
\item[a.] The set of isomorphism classes of $\mathcal{C}^1$ is
$HH^2_{-1}(A,A)$.
\item[b.] For $A^i\in \mathcal{C}^i$, the obstruction to its
extension to level $i+1$ lies in $HH^3_{-i-1}(A,A)$.
\item[c.] Furthermore, the set of isomorphism classes of such
extensions is an $HH^2_{-i-1}(A,A)$-homogeneous space.
\end{itemize}
\end{lemma}

For the proof, refer to \cite{BG96}.\\

We let $C^{\bullet}_{gr}(A,M)$ be the complex of graded morphisms from
the bar complex of $A$ to the $A$-bimodule $M$, so that in particular
the graded Hochschild cohomology groups appearing above are graded
cohomology modules of this complex.

Following Berger \cite{Ber01}, there is a complex
associated to an $N$-homogeneous algebra $A$, called its 
{\it $N$-Koszul} complex.  Its first few terms are

\[\begin{array}{l}

K_0=A\otimes A\\ K_1=A\otimes V \otimes A\\ K_2=A\otimes R\otimes A\\
K_3=A\otimes ((V\otimes R)\bigcap (R\otimes V))\otimes A
\end{array} \]

In general, the term $K_i$ has lowest degree terms $n(i)$, where
$n(2q)=Nq$ and $n(2q+1)=Nq+1$.
For details, including the definition of the differential, consult
\cite{Ber01}.

The algebra $A$ is $N$-Koszul iff this is a resolution of 
$A$ as an $A$-bimodule.
In other words, there is a
quasi-isomorphism from the Koszul complex to the bar complex:
\[\sigma:K_{\bullet}\rightarrow B_{\bullet}\]

It is injective since $K_{\bullet}$ is a minimal resolution, see \cite{Ber01}.
Note that this allows us to compute the graded Hochschild groups 
using maps from $K_{\bullet}\rightarrow A$; $\sigma$ induces a map
from $C^{\bullet}_{gr}(A,A)=Hom(B_{\bullet},A)$ to
$Hom(K_{\bullet},A)$, which is of course also a quasi-isomorphism.  
The differential we will be most interested in is the map 
$Hom(K_2,A)\rightarrow Hom(K_3,A)$.  Given a bimodule map 
$\alpha:K_2\rightarrow A$, 
its boundary $d\alpha$ is given on $(V\otimes R) \cap (R\otimes V)$ 
as $[\en,\alpha]$. \\

\begin{lemma}\label{hochnull} Suppose $A$ is  $N$-Koszul.  The graded
Hochschild cohomology groups $HH^i_j(A,M)$ vanish whenever $j < -n(i)$
and $M$ is a non-negatively graded $A$-bimodule.
\end{lemma}
This can easily be deduced from \cite{BG96} and \cite{BM03}; in fact,
in the spirit of \cite{BG96}, it might be reasonable to use this as
the definition of $N$-Koszulity.

We are now ready for the second part of the PBW-theorem.

\begin{proof}[Proof of Thm. \ref{NKThmPBW}, second part.]
The maps $\alpha_i:R\rightarrow V^{\otimes N-i}$ of vector spaces
naturally determine $A\otimes A^o$-maps from $K_2$ to $A$, which are
denoted by the same symbols.

To get the first level deformation, we need a class in
$HH^2_{-1}(A,A)$.  Condition (\ref{NKLikJac1}) from Theorem 
\ref{NKThmPBW} states
that the image of 
$[\en, \alpha_1]$ is contained in $R$ which is
equivalent to $d\alpha_1=0$ in the complex $Hom(K_{\bullet},A)$.  It thus
determines a class $[\alpha_1]$ in $HH^2_{-1}(A,A)$.  We want to choose
a representative for this class, $-f_1$, in the Hochschild complex,
such that $-f_1\circ\sigma=\alpha_1$ as cochains.  Here, and later, it
will always be possible to choose maps with the given properties
because the map $\sigma$ is determined by an (injective) map of vector
spaces.  So let $-\tilde{f_1}$
be any representative; then $-\tilde{f_1}\circ\sigma-\alpha_1$ is
cohomologically trivial, whence there exists $\omega$ in 
$Hom(K_1,A)_{-1}$ with $-\tilde{f_1}\circ\sigma-\alpha_1=d\omega$.  Let
$\omega'$ be any Hochschild cochain such that
$\omega'\circ\sigma=\omega$, and let $f_1=\tilde{f_1}+d\omega'$.  Then
\[-f_1\circ\sigma=-\tilde{f_1}\circ\sigma+d\omega=\alpha_1. \]

The next steps of the deformation are constructed inductively.  So
suppose the deformation has been constructed up to level $i$, meaning
that we have elements $f_1,...,f_i$ in $C^2_{gr}(A,A)$ such that, for
$j \leq i$, $-f_j\circ\sigma=\alpha_j$ and 
\[-df_j=f_{j-1}\circ [\en, f_1]\]
(analogous to the properties of the $\alpha$'s
from Theorem \ref{NKThmPBW}).  Here $f_0=0$.

We construct $f_{i+1}$ with these properties as follows: Let
\[-\tilde{f}_{i+1}\in C^2_{gr}(A,A)_{-i-1}\]
be any element such that $-\tilde{f}_{i+1}\circ\sigma=\alpha_{i+1}$.
Define 
\[ \gamma=-d\tilde{f}_{i+1}-f_i\circ [\en,f_1]. \]  
Then $\gamma\circ\sigma=0$, and since $\sigma$ is a
quasi-isomorphism, $\gamma$ is cohomologically trivial.  So there
exists $\mu$ in $C^2_{gr}(A,A)_{-i-1}$ with $d\mu=\gamma$.  Choose also
$\mu'$ (easy diagram chase) such that $\mu'\circ\sigma=\mu\circ\sigma$ and $d\mu'=0$.
Finally, set $f_{i+1}=\tilde{f}_{i+1}+\mu-\mu'$.  Then
\[-df_{i+1}=-d(\tilde{f}_{i+1}-\mu+\mu')=\gamma +f_i\circ 
[\en,f_1] -\gamma+0=f_i\circ [\en,f_1]\]
and
\[-f_{i+1}\circ\sigma=-\tilde{f}_{i+1}-\mu\circ\sigma+\mu'\circ\sigma=\alpha_{i+1}.\]

Compare this to the obstruction 
\[ \sum_{j=1}^{i}f_j\circ [\en,f_{i-j+1}] \]
for lifting the structure, coming from the
general theory of deformations.  All the terms with $j<i$ are cohomologically
trivial, since the composition with the quasi-isomorphism $\sigma$ is
zero for degree reasons.  The only remaining term is covered by
$-df_{i+1}$, as shown above, so the deformation can be extended to
level $i+1$.\\

This gives the deformation up to level $N+1$.  From there on, we use
Lemma \ref{hochnull}, which states that the cohomology groups
$HH^i_j(A,A)$ vanish for $j < -n(i)$. By Lemma \ref{HHbet}, 
for $i=3$ this says that the space
where the obstruction to extending the deformation further than level
$N+1$ lives is zero. For $i=2$ it says that the space of choices for
this extension is zero.\\

So we have a collection of $f_i$ such that the rule
\[a\times b=ab-\sum_{i=1}^{\infty}f_i(a,b)t^i\]
defines an associative multiplication on an algebra fibred over
$k[t]$ (this works since only finitely many terms are nonzero).  Let
$\tilde{U}$ be the fibre over $t=1$.  Since it is
constructed as a deformation, we know that $gr\,\tilde{U}=A$.  We have a
canonical map $\tilde{\phi}:V\rightarrow \tilde{U}$ which extends to a
map $\phi:U\rightarrow \tilde{U}$ (by definition of the $f_i$).  We
thus have a sequence of maps
\[A\stackrel{p}{\surj}\gr \,U \stackrel{gr \,\phi}{\rightarrow}\gr\,
\tilde{U}\stackrel{\iso}{\rightarrow}A\]
and the composition is the identity on $A$.  Thus $p$ is an
isomorphism.

\end{proof}

\eks It is easy to construct examples where the theorem
fails for non-Koszul algebras.  For instance, let $A$ be the
commutative algebra

\[A=k[x,y,z]/(y^2-xz,xy,z^2).\]

This algebra is not Koszul; in fact $HH_{-4}^3(A,A)\neq 0$.
It was considered in \cite{HRW98}.  Define $\alpha_1$ by $z^2\mapsto -x$, all
other relations (including the commutators) being mapped to zero.
Then $[\en,\alpha_1]$ is identically zero,
so the conditions of the theorem are trivially fulfilled.  But the
algebra $U=k[x,y,z]/(y^2-xz,xy,z^2+x)$ is finite dimensional as a
vector space, so $\gr \, U$ cannot be isomorphic to the algebra $A$,
which has Krull dimension 1.
\eksfin

\section{The $\ainf$-structure on the Yoneda algebra}

For any augmented algebra $A$ the Yoneda algebra $Ext_A^*(k,k)$ 
has a natural $\ainf$-structure, see \cite{KeMin}.
Keller computed this structure when $A$ is an 
$N$-Koszul algebra. To state this we let
$W=V^*$ and $S$ be the perpendicular space of $R$,
i.e. the kernel of $W^{\te N} \pil R^*$. 
Write $A^{!} = T(W)/(S)$ and let $B$ be the algebra with 
\[ B_{2p} = (A^!)_{Np}, \quad B_{2p+1} = (A^!)_{Np+1}.\]
Consider $k$ as a right $A$-module.
Then  $B$ is the Yoneda algebra $Ext_A^*(k,k)$ of $A$ as follows
from the $N$-Koszul complex for $A$. (If $k$ is considered as a left $A$-module
the Yoneda algebra will be the opposite algebra of $B$.)
Keller shows, \cite{KeNK}, that $m_2$ is ordinary
multiplication if some argument is even (in the $B$-grading).
When all arguments are odd then
$m_N(b_1, \ldots, b_N)$ is the ordinary product $b_1\cdots b_N$ in $B$.
In all other cases the $m_p$ vanish.

The following states that an augmented PBW-deformation $U$ of $A$ corresponds 
to giving an $\ainf$-structure on $B$ with certain explicitly
given properties.

\begin{theorem} \label{AinfTheEkv}
There is a one-to-one correspondence between PBW-defor\-mations of $A$ 
(with $\alpha_N = 0$) and $\ainf$-structures on $B$ with the following
properties: 
i) $m_n$ vanishes for $n > N$, 
ii) $m_N$
is the ordinary product of $N$ elements when all arguments are odd,
iii) $m_n$ vanishes for $3\leq n \leq N$ if some argument is even, and
iv) $m_2$ is ordinary multiplication when some argument is even.
\end{theorem}

\rem We see that the products in $B$ which are non-zero in the
non-deformed case are unchanged in the deformed case. The only 
products that we deform are the $m_p$ of odd arguments for
$p \leq N-1$ and $m_1$ of even arguments.
\remfin

 The above $\ainf$-structure on $B$ corresponds
to the natural $\ainf$-structure on the Yoneda algebra of $U$ as the 
following shows.
Recall that the bar construction $BU = \oplus_{n \geq 0} \overline{U}^{\te n}$,
where $\overline{U}$ is the kernel of $U \pil k$, comes
with a natural coalgebra structure, see \cite{Pr}.

\begin{theorem} \label{AinfTheYon} Let $C$ be the graded dual of $B$.
There is a quasi-isomorphism $C \pil BU$ of $\ainf$-coalgebras.
\end{theorem}

\rem It might be conceivable that there are other $\ainf$-structures on $B$
 corresponding to a given PBW-deformation, since $B$ is not in general
minimal. But the above Theorem \ref{AinfTheEkv} 
shows that there is a unique one with the above requirements.
\remfin

Let $\bev = \oplus B_{2p}$ be the {\it even} part of $B$, and 
$\bod = \oplus B_{2p+1}$ be the {\it odd} part of $B$.
We shall reserve $m_1$ to denote only the differential $\bod \pil \bev$.
The differential $\bev  \pil \bod$ shall be denoted by $d$. Also, if
some argument of $m_2$ is even, so $m_2$ is ordinary multiplication we
simply write $\cdot$ for this multiplication or even drop it so
we have the ordinary concatenation notation.

With the requirement that $m_p$ vanishes when $p \geq 3$ and some argument
is even, the axioms to be checked for $B$ to be an $\ainf$-algebra,
\cite[3.1]{KeInf} reduce to the following two types, where $p \geq 0$.
When all arguments are odd ($m_0 = 0$), 
\begin{equation}
 d\circ m_{p+1} = (-1)^{p}(\en\cdot m_p - m_p \cdot \en).
\label{AinfLigAks1}
\end{equation}
Let $u$ denote an even element and $a_i$ odd elements.
\begin{eqnarray}
\label{AinfLigAks2}
& & m_p(\ldots, a_{i-1},ua_i, \ldots)  \\
& = & m_p(\ldots,a_{i-1}u,a_i, \ldots)
 + (-1)^{p+1} m_{p+1}(\ldots,a_{i-1}, d(u), a_i, \ldots ). \notag
\end{eqnarray}
(If $i = 1$ then in the middle term $u$ gives ordinary multiplication to 
the left "outside" $m_p$, and similarly for the first term when $i = p+1$.)

\rem An interesting thing here is that one realizes that this may
be abstracted to give the following construction of a class of 
$\ainf$-algebras. Let $\bev$ be a positively even graded algebra, and
$\bod$ a module graded in positive odd degrees, with a derivation 
$d$ and maps $m_p$ of degree $2-p$
\[ \bev \mto{d} \bod, \quad
 (\bod)^{\te p} \mto{m_p} \bev \]
such that (\ref{AinfLigAks1}) and (\ref{AinfLigAks2})
hold. Then $B = \bev \oplus \bod$ becomes an $\ainf$-algebra.
\remfin

\rem
Note that when $p=1$ both (\ref{AinfLigAks1}) and (\ref{AinfLigAks2})
are instances of the derivation rule $\delta(ab) = \delta(a)b + 
(-1)^{|a|}a \delta(b)$. Hence they generalize this equation in different
directions.
\remfin

\medskip

We shall now show how the $\ainf$-structure on $B$ determines a 
PBW-deformation.

\begin{proof}[Proof of Theorem \ref{AinfTheEkv}, the easy direction.]
The $\ainf$-structure on $B$ provides maps
\begin{equation}
W^{\te p} = (B_1)^{\te p} \mto{m_p} B_2 = W^{\te N}/S. \label{AinfLigMp}
\end{equation}
We define $\alpha_{p} : R \pil V^{\te N-p}$ to be the sign-adjusted dual
given by
\begin{equation}
\alpha_{p} = \sigma(p) \cdot (m_{N-p})^*
\end{equation}
where $\sigma(p)$ is a sign determined by $\sigma(p-2) = - \sigma(p)$, 
$\sigma(N)= 1$ and $\sigma(N-1) = (-1)^{N-1}$.
The equation (\ref{AinfLigAks1}) for $p < N-1$ translates to the
condition (\ref{NKLikJac2}) of Theorem \ref{NKThmPBW}. 
The equation (\ref{AinfLigAks1}) for $p=N-1$
shows that the map
\begin{equation} 
W^{\te N} \pil W^{\te N}/S \mto{d} W^{\te N+1}/(S\te W + W \te S) 
\label{AinfLigDif}
\end{equation}
is the dual of the map
\[ (V \te R) \cap (R \te V) \mto{[\en,\alpha_1]} R \sus V^{\te N} \]
showing condition (\ref{NKLikJac1}) of Theorem \ref{NKThmPBW}.
\end{proof}

Considerably more work it is to show the other direction of 
Theorem \ref{AinfTheEkv}.
We shall do this through several steps in the next section. Here we shall
explain how, given the $\alpha_i$'s, to define the maps $d$ and $m_p$.

\medskip

\noindent{\bf Definition of $d$.} First we define 
\[ W^{\te N} = B_2 \mto{d} B_3 \] as the dual
of 
\[ (V\te R) \cap (R\te V) \mto{[\en,\alpha_1]} V^{\te N}. \]
This may be extended to a derivation from $T(W)^{ev} = \oplus_i T(W)_{Ni}$
to $\bod$. We shall 
show that this descends to a map $\bev \pil \bod$.

\medskip

\noindent{\bf Definition of $m_p$.} When all arguments are linear
define $m_p$ to be the sign-adjusted dual
\[ m_p = \sigma(N-p)(\alpha_{N-p})^*. \]
Supposing $m_{p+1}$ has been defined, we may
then define 
\[ (T(W)^{od})^{\te p} \mto{m_p} \bev \]
inductively by equation (\ref{AinfLigAks2})
on the smallest integer $i$ such that $a_{i+1}, \ldots, a_p$ are all linear.
Then let $ua_i$ be a factorization with $a_i$ linear.
We shall show that this descends to a map $(\bod)^{\te p} \pil \bev$.

\medskip
In the next section we shall prove that these definitions provide
$B$ with an $\ainf$-structure. Assuming this we here provide the proof
of Theorem \ref{AinfTheYon}.

\begin{proof}[Proof Theorem \ref{AinfTheYon}]
The map $C \mto{\tau} U$ which is the projection on $C_1$
composed with the inclusion $C_1 = V \inpil U$
is a (generalized) twisting cochain, i.e. we have
\[ \sum_{p=1}^\infty \mu_U \circ \tau^{\te p} \circ (m_p)^* = 0 \]
where $\mu_U$ is the iterated multiplication on $U$.
In fact this reduces simply to the equation (note that $\tau$ has degree one)
\[ r + \sum_i \alpha_i(r) = 0.\]
Hence we get an $\ainf$-coalgebra morphism $C \pil BU$.
As in \cite{KeMin}, see also \cite{KeNK}, in order to show that this is
a quasi-isomorphism, it will be sufficient to show that the twisted
tensor product
\[ C\te_{\tau} U = C_0 \te U \mlto{d_\tau} C_1 \te U \mlto{d_\tau} \cdots \]
is a resolution of $k$, where the differential $d_{\tau}$ is given by
\[ \sum_{p=1}^\infty (\one_C \te \mu_U)(\one_C \te \tau^{\te {p-1}} \te 
\one_U) ((m_p)^* \te \one_U). \]
But the complex $C \te_{\tau} U$ comes with a natural filtration
$(C \te_{\tau} U)_{\leq p}$ given by
\begin{eqnarray*} 
 & C_0 \te U_{\leq p} &\vpil C_1 \te U_{\leq p-1} \vpil C_2 \te U_{\leq p-N}
\vpil C_3 \te U_{\leq p-N-1} \\
 \vpil & C_4 \te U_{\leq p-2N} & \vpil \cdots .
\end{eqnarray*}
Now the "non-deformed" graded complex $C\te_{\tau_0} A$,
where $\tau_0$ is the twisting cochain when $U$ is specialized to $A$,
is the Koszul
complex for the $N$-Koszul algebra $A$ and is thus a resolution of $k$.
It is easily worked out that, for each $p$, the quotient 
\[ (C\te_{\tau} U)_{\leq p} / (C \te_{\tau} U)_{\leq p-1} \]
identifies as the strand of $C \te_{\tau_0} A$ in degree $p$, and so is acyclic
for $p \geq 1$. The upshot is that $(C\te_{\tau} U)_{\leq p}$ is 
quasi-isomorphic to $k$ for all $p \geq 0$. Since $C \te_{\tau} U$ is 
the colimit of the $(C \te_{\tau} U)_{\leq p}$ and this colimit commutes with
homology, we get that $C\te_{\tau} U$ is a resolution of $k$.
\end{proof}

\section{Proof of Theorem \ref{AinfTheEkv}}

Assuming the definitions of $d$ and $m_p$, we now proceed to show that
the maps $d$ and $m_p$ descend, and that equations (\ref{AinfLigAks1}) and
(\ref{AinfLigAks2}) hold. Our first lemma will in particular be of help in 
showing the descent of the maps $d$ and $m_p$. Our second lemma
shows the descent of $d$. Then we proceed through a sequence of lemmas
\ref{ProLemAks2} to \ref{ProLemDesmp}
concerning the $m_p$. The argument is here inductive in the sense
that each time we prove a lemma concerning $m_p$ we assume that {\it all} 
the lemmas have been proven for $m_r$ when $r \geq p+1$.

\begin{lemma}  \label{ProLemWS} For $m = 2,3, \ldots, N-1$
\begin{equation}
W^{\te m-1} \te S \te W \sus W^{\te m} \te S + \cap_{i=1}^m (W^{\te m-i}
\te S \te W^{\te i}). \label{ProLigWS}
\end{equation}
\end{lemma}

\begin{proof} From \cite{Ber01} we know that 
\begin{equation} 
  (V^{\te m}\te R) \cap (R \te V^{\te m} + V \te R \te V^{\te m-1} + 
\cdots + V^{\te m-1} \te R \te V) \label{ProLigVR1}
\end{equation}
is equal to
\begin{eqnarray}
 & & (V^{\te m} \te R) \cap (R \te V^{\te m}) \label{ProLigVR21} \\
&+&  (V^{\te m} \te R) \cap 
( V \te R \te V^{\te m-1} + \cdots +  V^{\te m-1} \te R \te V) 
\label{ProLigVR22}
\end{eqnarray}
and also 
\[ (V^{\te m} \te R) \cap (R \te V^{\te m}) \sus V^{\te m-1} \te R \te V. \]
When $m=2$ this gives that (\ref{ProLigVR1}) is included in $V \te R \te V$.
When $m > 2$ we may by induction assume that (\ref{ProLigVR22}) is 
included in $V \te V^{\te m-2} \te R \te V$. Whence we get
\begin{eqnarray}
\label{ProLigVR}
 & & (V^{\te m} \te R) \cap (R \te V^{\te m} + V \te R \te V^{\te m-1} + \cdots
+ V^{\te m-1} \te R \te V) \\ 
& \sus &  V^{\te m-1} \te R \te V. \notag
\end{eqnarray}

In general when $V$ and $W$ are dual vector spaces and $I$ and $I^\bot$
and $J$ and $J^\bot$ are orthogonal complements for the pairing, then
$I^\bot + J^\bot$ will be the orthogonal complement of $I \cap J$. 
Taking the orthogonal
complement of (\ref{ProLigVR}), gives the statement of the lemma.
\end{proof}

\begin{lemma}
The map $T(W)^{ev} \mto{d} \bod$ descends to a map $\bev \pil \bod$.
\end{lemma}

\begin{proof}
Note first that $d(S) = 0$. This is because $W^{\te N} \pil B_3$
is dual to $(V\te R) \cap (R \te V) \mto{[\en,\alpha_1]} 
V^{\te N}$ which factors
through $R \sus V^{\te N}$.
Since $d$ is a derivation, it will be enough to show that 
$d(W^{\te a} \te S \te W^{\te b})$ becomes zero when $a+b = N$.
Note that this vanishes by the derivation property if $a$ or $b$ is zero.

Suppose $b = 1$. By tensoring Eq. (\ref{ProLigWS}) to the left with
$W$ when $m = N-1$ we get
\begin{equation}
W^{\te N-1} \te S \te W \sus W^{\te N} \te S + \cap_{i=1}^{N-1} (W^{\te N-i}
\te S \te W^{\te i}). \label{ProLigWNm1}
\end{equation}
It will therefore be enough to show that $d$ applied to the last term in
(\ref{ProLigWNm1}) is zero.

Before proceeding note the following.
Since $W^{\te N} \mto{d} B_3$ is the dual of
\[ (V\te R) \cap (R \te V) \mto{ [\en,\alpha_1]} V^{\te N} \]
we have the equation (when arguments are in $W$)
\[ d m_N = (-1)^{N-1}( \one \cdot m_{N-1} - m_{N-1}\cdot \one ). \]
Summing up, we get (with arguments in $W$)
\begin{equation} 
\sum_{i=0}^{N} \ent{i} \cdot d m_N \cdot \ent{N-i} 
= (-1)^{N-1} (m_{N-1} \cdot \one^{\te N+1} - \one^{\te N+1} \cdot m_{N-1}). 
\label{ProLigDmN}
\end{equation}
Applying this to an argument in the last term of (\ref{ProLigWNm1})
the right side of (\ref{ProLigDmN}) becomes
zero in $B_5 = (A^!)_{2N+1}$ and the left side reduces to 
\[ dm_N \cdot m_N + m_N \cdot d m_N. \]
But this is equal to $d(m_N \te m_N)$ by the derivation property, and hence
$d$ vanishes on  the last term of (\ref{ProLigWNm1}).

When $b = 2$ then tensoring Eq.(\ref{ProLigWS}) to the right with $W$ when
$m = N-1$ we get
\begin{equation*}
 W^{\te N-2} \te S \te W^{\te 2} 
\sus  W^{\te N-1} \te S \te W + S \te W^{\te N}
\end{equation*}
and hence $d(W^{\te N-2} \te S \te W^{\te 2})$ is zero.
In this way we may continue tensoring (\ref{ProLigWS}) to the right with
powers of $W$ and get inductively that $d(V^{\te a} \te S \te W^{\te b})$
is zero.
\end{proof}

\begin{lemma}\label{ProLemAks2}
 Equation (\ref{AinfLigAks2}) holds when $i \leq p$.
\end{lemma}

\begin{proof} For the induction start note that it holds when $p = N$.
Now if the $a_j$ are linear for $j \geq i$, 
Eq. (\ref{AinfLigAks2}) is just the
definition.  So let $j \geq i$ be the least index such that $a_{j+1}, \ldots,
a_p$ are all linear, but $a_j$ is not linear.

Suppose first that $j \geq i+2$ and let $a_j = vl$ where $v$ is even
(i.e.
in $W^{\te Nq}$ for some $q$), and $l$ is linear. The first term
in Eq. (\ref{AinfLigAks2}) is then by definition equal to 
\begin{eqnarray}\label{ProLigUV1}
& &  m_p( \ldots, a_{i-1}, ua_i, \ldots, a_{j-1}v,l, \ldots) \\
& +  (-1)^{p+1} 
&  m_{p+1}(\ldots, a_{i-1}, ua_i, \ldots, a_{j-1},d(v),l, \ldots)
\notag 
\end{eqnarray}

The second term in Eq.(\ref{AinfLigAks2}) is by definition
\begin{eqnarray} \label{ProLigUV2}
& & m_p(\ldots, a_{i-1}u, a_i, \ldots , a_{j-1}v,l, \ldots) \\ 
& + (-1)^{p+1}& 
m_{p+1}(\ldots, a_{i-1}u, a_i, \ldots, a_{j-1}, d(v), l, \ldots ) \notag
\end{eqnarray}
and the third term in Eq. (\ref{AinfLigAks2}) is 
\begin{eqnarray} \label{ProLigUV3}
 & 
(-1)^{p+1} & m_{p+1}(\ldots, a_{i-1}, d(u), a_i, \ldots, a_{j-1}v,l, \ldots) \\
& - & m_{p+2}(\ldots,a_{i-1}, d(u), a_i, \ldots, a_{j-1}, d(v), l, \ldots ).
\notag
\end{eqnarray}
By ascending induction on $j$ (we'll take care of the initial cases $j = i, i+1$
shortly) and the induction assumptions, 
we may assume that (\ref{ProLigUV1}) is 
\begin{eqnarray*}
& &  m_p(\ldots, a_{i-1}u, a_i, \ldots, a_{j-1}v,l, \ldots) \\ 
& + (-1)^{p+1} & m_{p+1}(\ldots, a_{i-1}, d(u), a_i, \ldots, a_{j-1}v, l, \ldots)\\
& + (-1)^{p+1} & m_{p+1}(\ldots, a_{i-1}u, a_i, \ldots, a_{j-1}, d(v), l, \ldots ) \\
& - & m_{p+2}(\ldots, a_{i-1}, d(u), a_i, \ldots, a_{j-1}, d(v), l, \ldots)
\end{eqnarray*}
which is the sum of (\ref{ProLigUV2}) and (\ref{ProLigUV3}).

When $j = i+1$ the same type of argument goes through, assuming it is 
shown for $j = i$, and when $j = i$, an analogous argument goes through,
this time using $d(uv) = d(u)v + ud(v)$.
\end{proof}

\begin{lemma} Equation (\ref{AinfLigAks1}) holds. \label{ProLemAks1}
\end{lemma}

\begin{proof}
If all arguments are linear, this is simply the dual of the conditions
(\ref{NKLikJac2}) of Theorem \ref{NKThmPBW}. An odd
argument which is not linear can be written (as a sum of)
$m_N(a_1, \ldots, a_N)\cdot a_{N+1}$ where the $a_i$ are odd elements. Suppose
the $l+1$'th argument of Eq. (\ref{AinfLigAks1}) 
is of this form, where $l \leq p$.
We must then show (note that we write $m_N \te \one$ because we have 
not yet shown that
$m_p$ descends, see Lemma \ref{ProLemDesmp})
\begin{eqnarray}\label{ProLigDMp}
& & d m_{p+1}(\one^{\te l} \te m_N\cdot \one \te \one^{\te p-l}) \\
& =  (-1)^{p}&
(\one \cdot m_p(\one^{\te l-1} \te (m_N \te \one) \te \one^{\te p-l})\notag \\
& - &  m_p(\one^{\te l} \te 
(m_N \te \one) \te \one^{\te p-1-l})\cdot \one) \notag 
\end{eqnarray}
Considering the first term, we have by using Eq. (\ref{AinfLigAks2}),
successively ``shifting the $u$ term all the way to the left'', that
$m_{p+1}(\one^{\te l} \te m_N\cdot \one \te \one^{\te p-l})$ is
\begin{equation}
m_N \cdot m_{p+1}(\one^{\te p+1}) + (-1)^{p+2} \sum_{i = 0}^l
m_{p+2}(\one^{\te i} \te dm_N \te \one^{\te p+1-i})
\label{ProLigMpshift}
\end{equation}
Applying $d$ to this, the first term in (\ref{ProLigDMp}) is
\begin{eqnarray}
& & dm_N \cdot m_{p+1}(\one^{\te p+1}) + m_N \cdot dm_{p+1}(\one^{\te p+1}) 
\label{ProLigDMps1}\\
& -  &\sum_{i=1}^l \one \cdot m_{p+1}(\one^{\te i-1} \te dm_N 
\te \one^{\te p-i}) 
\label{ProLigDMps3} \\
& - & dm_N \cdot m_{p+1}(\one^{\te p+1}) \label{ProLigDMps4} \\
& +  &\sum_{i=0}^l  m_{p+1}(\one^{\te i} \te dm_N \te \one^{\te p-i})
\cdot \one \label{ProLigDMps2}
\end{eqnarray}
The first term in (\ref{ProLigDMps1}) and the term in (\ref{ProLigDMps4})
cancel. By induction on the total degree of the argument, the second
term in (\ref{ProLigDMps1}) is equal to
\begin{equation} (-1)^{p} m_N \cdot (\one \cdot 
m_p(\one^{\te p}) - m_p(\one^{\te p})\cdot \one). \label{ProLigDMps5}
\end{equation}
But then again by using Eq. (\ref{AinfLigAks2}) and 
``shifting successively to the right'', we see that the last two terms
in (\ref{ProLigDMp}) become equal to the sum of 
(\ref{ProLigDMps5}) and the terms in (\ref{ProLigDMps2}) and 
(\ref{ProLigDMps3}).
\end{proof}

\begin{lemma}  \label{ProLemMpN} 
\begin{equation}
 \sum_{i=0}^p m_{p+1}(\one^{\te i} \te dm_N \te \one^{\te p-i})
= (-1)^{p}(m_N \cdot m_p - m_p \cdot m_N) \label{ProLigMpN}
\end{equation}
\end{lemma}

\begin{proof}
According to Lemma \ref{ProLemAks1} we can assume that the equation
\[ d m_{r+1} = (-1)^{r} (\one \cdot m_r - m_r \cdot \one) \]
holds for $r \geq p$. Assuming $q \geq p+1$ and inserting this in the sum
\[ \sum_{i = 0}^{q-1}m_q(\one^{\te i} \te dm_{r+1} \te \one^{\te q-1-i}) \]
we get, observing that many terms will cancel because of Eq. 
(\ref{AinfLigAks2}), that this is equal to 
\[(-1)^{r+q} \sum_{i=0}^q m_{q+1}(\one^{\te i} \te dm_r \te \one^{\te q-i} )
+ (-1)^r(m_q m_r - m_r m_q). \]
Note that the first term here is of the same type which we just
expanded, so we may proceed inductively.
We now start with the left hand side of (\ref{ProLigMpN})
and apply this inductively. In the end this becomes
\begin{equation}
\sum_{a\geq p+1, a+b = N+p} (-1)^{\eta(a,b)} (m_am_b - m_b m_a)
\label{ProLigMab}
\end{equation}
where $\eta(a,b) \equiv \eta(b,a) (\text{mod } 2)$. Hence the terms in 
(\ref{ProLigMab}) cancel, except the term
\[(-1)^{p} (m_N \cdot m_p - m_p \cdot m_N). \]
This proves the lemma.
\end{proof}

\begin{lemma} Equation (\ref{AinfLigAks2}) holds when $i= p+1$.
\end{lemma}

\begin{proof} By the fact that Eq.(\ref{AinfLigAks2}) holds for $i \leq p$, 
and successively
``shifting even factors all the way to the right outside of $m_p$'', we
are reduced to prove the $i = p+1$ case when all the $a_i$  are linear.
Also $u$ may be written as $m_N(a_{p+1}, \ldots, a_{N+p})$ 
where these arguments
are odd. The second term in (\ref{AinfLigAks2}), $m_p(\ldots, a_pu)$ will then,
by successively ``shifting even terms all the way to the left and then
outside of  $m_p$'',
be equal to
\[a_1 \cdots a_N m_p(a_{N+1},\ldots, a_{N+p}) + (-1)^{p+1}
\sum_{i=1}^{p+1} m_{p+1}(\ldots, d(a_i \cdots a_{i+N-1}), \ldots). \]
By the previous Lemma \ref{ProLemMpN}, this is equal to
$m_p(a_1 \cdots a_p) a_{p+1} \cdots a_N$.
\end{proof}

\begin{lemma} \label{ProLemDesmp}
The map 
\[ (T(W)^{od})^{\te p} \mto{m_p} \bev\]
descends to a map $(\bod)^{\te p} \pil \bev$.
\end{lemma}

\begin{proof}
We must show that $m_p$ becomes zero if one of its $p$ arguments
have the form $W^{\te a} \te S \te W^{\te b}$. 
Using Eq.(\ref{AinfLigAks2}) and
shifting even factors to the left or right outside of $m_p$, this is
easily reduced to the case when all arguments are linear, except one
which has the form $W^{\te a} \te S \te W^{\te b}$ where $a+b$ is $N+1$
or $1$. 
Suppose $a=0$ (so $b=1$) and that the non-linear argument is in position
$p$, the last one.
We must show that $m_p$ vanishes on $W^{\te p-1} \te (S \te W)$. By
Lemma \ref{ProLemWS} this is contained in
\[ W^{\te p} \te S + \cap_{i=1}^p (W^{\te p-i} \te S \te W^{\te i}) \]
and hence $W^{\te p-1} \te (S \te W)$ is contained in 
\begin{equation} 
(W^{\te p-1} \te S \te W) \cap (W^{\te p} \te S) + \cap_{i = 1}^p 
(W^{\te p-i} \te S \te W^{\te i}). \label{ProLigWpS}
\end{equation}
Now if the argument of $m_p$ is contained in the left term of 
(\ref{ProLigWpS}), then $m_p$ vanishes on it by using Eq.
(\ref{AinfLigAks2})
and ``shifting $S$ to the right outside of $m_p$''. If the argument of  
$m_p$ is contained in the right term of (\ref{ProLigWpS}) we may
use Eq. (\ref{AinfLigAks2}) again and shift the argument of $m_p$ 
successively to 
the left. Note that all terms with $m_{p+1}$ will vanish since $d(S)$ is
zero. In the second last step we are left with an argument in 
$(S \te W) \te W^{\te p-1}$ and in the last step we ``shift the $S$-part
to the left outside of $m_p$'', and so $m_p$ vanishes on this argument.

Assume now that it is the argument in position $p-1$ of $m_p$ that
is contained in $S\te W$. By Lemma \ref{ProLemWS}, 
$W^{\te p-2} \te (S \te W) \te W$ is equal to 
\begin{eqnarray*}
 & & (W^{\te p-2} \te (S\te W) \te W) \cap (W^{\te p-1} \te (S \te W)) \\
& + &
\cap_{i=1}^{p-1} (W^{\te p-1-i} \te (S \te W) \te W^{\te i}). 
\end{eqnarray*}
By the former case $m_p$ vanishes on the first term and by the same kind
of shift argument as above, it also vanishes on the second term. In this
way we may continue, and this settles the $a=0$ case.

The case $a=1$ and $b=0$ may be reduced to the case just treated, by shifting
the $S$-part one step to the right, using Eq. (\ref{AinfLigAks2}).

So suppose now $2 \leq a \leq N-1$ and $a+b = N+1$. Applying Lemma 
\ref{ProLemWS} with $a= m-1$ and tensoring Eq.(\ref{ProLigWS})
to the right with $W^{\te b-1}$, we derive
\[ W^{\te a} \te S \te W^{\te b} \sus W^{\te a+1} \te S \te W^{\te b-1}
+ S \te W^{\te N+1}. \]
So by induction we get that $m_p$ vanishes when this is one of the arguments.
\end{proof}

\section{PBW-deformations of $T(V)/(\wedge^N V)$.}

In \cite{Ber01}, R. Berger showed that the algebra $T(V)/(\wedge^N V)$, where
$\wedge^N V$ is naturally embedded in  $V^{\te N}$ by
\[ x_1 \wedge \cdots \wedge x_N \mapsto \sum_{\sigma \in S_N} (-1)^{sgn
(\sigma)} x_{\sigma(1)} \te \cdots \te x_{\sigma(N)} \]
is an N-Koszul algebra (char $k = 0$). Theorem \ref{NKThmPBW} 
therefore applies
and we shall use this to compute the possible $\alpha_i$'s for this algebra
such that the associated algebra $U$ is of PBW-type.

When $N=2$ (and $\alpha_2 = 0$) the classical PBW-theorem says that $U$ is 
of PBW-type
iff $\alpha_1 : \wedge^2 V \pil V$ is a Lie bracket. In the general case
when $N$ is even
the result involves a Lie algebra together with symplectic forms
of successive even degrees and the Lie algebra is related to each form
by an equation which can be viewed as generalizing the Jacobi identity.

The case $N$ odd is distinct and simply involves choosing
freely a linear form and 
successive symplectic forms of even degree. 
We shall now state the main theorem, but first
we explain some notation.

\medskip
\noindent {\bf Notation.} 1. If we have a map 
\begin{equation} \psi : \wedge^p V \pil U \label{SetNotP1}
\end{equation}
the domain may be considered as a quotient space of $V^{\te p}$.
Thus we get a map which we denote by the same symbol
\begin{equation} 
\psi : V^{\te p} \pil U. \label{SetNotP2}
\end{equation}
On the other hand given map (\ref{SetNotP2}) then by restricting it to
the subspace $\wedge^p V$ of the domain we get a map (\ref{SetNotP1}).
We shall be free to switch back and forth like this.
Thus if $\psi^\prime : \wedge^q V \pil U^\prime$ is another
map we may form $\psi \te \psi^\prime$ which we may consider as 
defined on $\wedge^{p+q} V$.
If we go from (\ref{SetNotP1}) to (\ref{SetNotP2}) and back to   
(\ref{SetNotP1}) again the new map will differ from the old
by a constant, but only in a few cases will it be of any importance to
keep track of this.

\medskip
2.    There is also a natural map
\begin{equation} 
\Hom_k( \wedge^p V, \wedge^r V) \mto{T_a} 
   \Hom_k( \wedge^{a+p} V, \wedge^{a+r} V) \label{Se-2}
\end{equation}
given by (arrange $I,P$ and $A$ in ascending order)
\[ T_a(\phi)(v_I) = \sum_{I = A \cup P} (-1)^{sgn(A,P)} 
 v_A \wedge  \phi(v_P), \]
where we sum over all partitions $A \cup P$ of $I$ with $P$ of 
cardinality $p$.
The sign is the sign of the permutation we get by concatenating $P$ with $A$.
Given a map $\phi$ we shall often by abuse of notation denote $T_a(\phi)$
simply by $\phi$, like in (\ref{SeLaGL}) below.

\medskip 
3. Let $\ent{p}$ denote the identity on $V^{\te p}$.
For a linear map $L : V \te V \pil V$ we introduce the notation
\[\underline{\en^{2a}L^b} = \underset{\sum i_l = a,\sum j_l = b}{\sum}
  \ent{2i_1} \te  L^{\te j_1} \te \cdots \te \ent{2i_s} \te L^{\te j_s} \]
which is a map from $V^{\te  2a+2b }$ to $V^{\te  2a+b}$, and
\[ \underline{\smpm  \en^a L}  = 
\ent{a} \te L -   \ent{a-1} \te L \te \en +  \cdots + (-1)^a L \te \ent{a}. \]
Note that 
\[ [\en, \underline{\en^{2c} L}] 
= \underline{ \smpm \en^{2c+1} L}. \]
On some occasions we will also allow $L$ to be a linear map from 
$V^{\te s}$ to $V$ (notably in Proposition \ref{RePrSt}).

\medskip  
Now we shall consider PBW-deformations $U$ of the algebra $T(V)/(\wedge^N V)$
given by linear maps $\alpha_i : \wedge^N V \pil V^{\te N-i}$.

\begin{theorem} \label{Se-T1}
Suppose $N$ is odd and $\dim V \geq N+2$. Let $l : V \pil k$ be a linear
map and $\Phi_{2r} : \wedge^{2r} V \pil k$ arbitrary symplectic forms
for $2 \leq 2r < N$. Let
\begin{eqnarray}
\alpha_{2r} &=& \ent{N-2r} \te \Phi_{2r} \notag \\
\alpha_{2r+1} &=& \ent{N-2r-1} \te l \te \Phi_{2r} \label{SeLigAll}
\end{eqnarray}
Then these $\alpha_i$'s give an algebra $U$ of PBW-type. 
Conversely if $U$ is of
PBW-type there are $l$ and $\Phi_{2r}$ such that the $\alpha_i$ are of the form
above.
\end{theorem}

\begin{theorem} \label{Se-T2}
Suppose $N = 2n$ is even and $\dim V \geq N+2$. Let 
$L : \wedge^2 V \pil V$ be a Lie
bracket and $\Phi_{2r} : \wedge^{2r} V \pil k$ for $0 \leq 2r < N$ 
be symplectic forms (with $\Phi_0 = 1$) such that
the compositions 
\begin{equation} L \circ \Phi_{2r} \circ L : \wedge^{2r+3} V \pil V 
\label{SeLaGL} \end{equation} 
are zero. Let $\alpha_{2r}$ be
\begin{equation}
\sum_{i=0}^r \, \, \underline{\en^{2(n-2r+i)}L^{2r-2i} } \te \Phi_{2i} 
\label{Se-5}
\end{equation}
and let $\alpha_{2r+1}$ be
\begin{equation}
\sum_{i = 0}^r \, \, \underline{\en^{2(n-2r-1+i )}L^{2r+1-2i}} 
\te \Phi_{2i} + 
\ent{2n-2r-1} \te r\Phi_{2r} (\ent{2r-1}\te L)    \label{Se-6} 
\end{equation}
where $2r,2r+1$ assumes values from $1$ to $N-1$ and $\alpha_N = 0$. Then these
$\alpha_i$ give an algebra $U$ of PBW-type. Conversely if $U$ is of PBW-type
and $\alpha_N = 0$, there exists $L$ and $\Phi_{2r}$ such that the $\alpha_i$
are of the form above.

If $\alpha_N \neq 0$, the $U$ of PBW-type are obtained by also 
choosing a $\Phi_{2n} : \wedge^{2n} V \pil k$ such that the following 
composition is zero
\[\Phi_{2n} (\ent{2n-1} \te L) : \wedge^{2n+1} V \pil k.\]

\end{theorem}

\rem If the dimension of $V$ is $\leq N$ the $\alpha_i$ may be chosen 
arbitrarily.

\rem The case when $V$ has dimension $N+1$ seems quite different and 
would have needed a
separate treatment. So we do not consider it here.

\rem When $N = 4$ and $\alpha_4 = 0$ we have given a Lie bracket $L$ and 
a quadratic symplectic form $\Phi_2$ such that the composition 
$L \circ \Phi_2 \circ L$ is zero. 
Up to coordinate change a quadratic symplectic
form is given by its rank, an even number. For each such rank we thus get 
a class of Lie algebras. 
\remfin

Having stated the above results for the algebra $T(V)/(\wedge^N V)$, one
might wonder about what happens concerning PBW-deformations of the algebra
$T(V)/(S^N(V))$, the relations being the symmetric tensors
of degree $N$. We state the following without proof (the case $N=2$ is
well known for Clifford algebras).

\begin{theorem} Consider the algebra $T(V)/(S^N(V))$. To give linear maps
$S^N(V) \mto{\alpha_i} V^{\te N-i}$ such that $U$ is of PBW-type is
equivalent to giving arbitrary symmetric forms
$\Phi_r : S^r(V) \pil k$ for $r = 1, \ldots, N$.
(The $\alpha_i$'s are then given by suitable expressions in these
forms.)
\end{theorem}

\medskip

Before embarking on the proofs of the theorems above we shall need            
auxiliary results. 
We shall gradually move towards the proof of these 
theorems by investigating the conditions (\ref{NKLikJac1}) and
(\ref{NKLikJac2}) of Theorem \ref{NKThmPBW}. 
Concerning notation, a monomial $\wedge^p V$ 
will be denoted by $x_1 \wedge \cdots \wedge x_p$ while a monomial in
$V^{\te p}$ will be denoted simply by $x_1 \cdots x_p$.

\section{Results on the $\alpha$-mappings}

This section is mainly devoted to study linear maps
\[ \wedge^N V \mto{\alpha} V^{\te N-s}. \]
But first we need a result on the map $T_a$ given in (\ref{Se-2}).

\subsection{The linear map $T_a$.}

\begin{proposition} \label{RePrTa}
If $v=\dim_k V  \geq a+p+r$ then $T_a$ is injective.
In particular it is an isomorphism when $\dim_k V = a+p+r$.
\end{proposition}

\begin{proof}
Fix a basis $e_1, \ldots, e_v$ of $V$ and let $V = V^\prime \oplus 
\langle e \rangle$ where $e$ is one of the $e_i$'s. 
There is a decomposition $\wedge^p V = \wedge^p V^\prime 
\oplus \wedge^{p-1} V^\prime \wedge e$. For short denote
$\Hom_k(\wedge^s V^\prime, \wedge^t V^\prime)$ by $H(s,t)$.
There is then a decomposition of 
$\Hom_k(\wedge^p V, \wedge^r V)$ 
\[ H(p,r) \oplus H(p,r-1) \oplus H(p-1,r)
\oplus H(p-1,r-1) \]
and a corresponding decomposition for $\Hom_k(\wedge^{a+p} V, \wedge^{a+r} V)$.

Check that  $H(p,r-1)$ only maps to $H(a+p,a+r-1)$. By induction
this map is injective. Also  $H(p-1,r)$ only maps to 
$H(a+p-1,a+r)$ and so by induction is injective. The kernel of $T_a$
is therefore contained in $H(p,r) \oplus H(p-1,r-1)$. Take 
note that this maps to the subspace 
$H(a+p,a+r) \oplus H(a+p-1,a+r-1)$. 

For a subset $I$ of $\{1, \ldots, v\}$ denote by $e_I$ the product
$\wedge_{i \in I}e_i$. 
Now the maps $\phi_{P,R}$ in $\Hom_k(\wedge^p V \wedge^r V)$, sending
$e_P$ to $e_R$ and all other $e_{P^\prime}$ to zero,
constitute a basis for this space.
Suppose then $\sum \ga_{P,R} \phi_{P,R}$ is in the kernel of $T_a$.
Let $e$ be in $P\backslash R$. Then by what was said just above 
$\ga_{P,R} = 0$.
Similarly if we choose an $e$ in $R \backslash P$. 
Hence if $p \neq r$ we are done. We may therefore assume that $p = r$ and
that the kernel
consists of elements $\sum \ga_{P,P} \phi_{P,P}$.

So let $I(p)$ be the subspace of $\Hom_k(\wedge^p V, \wedge^p V)$
generated by $\phi_{P,P}$. Then $T_a$ restricts to a map 
$I(p) \pil I(p+a)$ which may be identified as follows.
Let $A$ be the commutative ring $k[e_1, \ldots, e_v]/(e_1^2, \ldots, e_v^2)$. 
We can identify $I(p)$ with $A_p$ via
$e_P \mapsto \phi_{P,P}$. If $\sigma_a$ is the $a$'th elementary 
symmetric polynomial in the variables $e_i$, then $T_a$ can 
be identified with the map 
\[ A_p \mto{\cdot \sigma_a} A_{a+p}. \]
Note that $a! \sigma_a = (\sigma_1)^a$ in  $A$. Also 
$\sigma_b \cdot \sigma_a = \binom{a+b}{a} \sigma_{a+b}$. 
To show that $\sigma_a$ is injective
it is therefore enough to show that $\sigma_{a+b}$ is injective
where the dimension of $V$ is $a+b+2p$.
Replacing $a$ with $a+b$ and $\sigma_a$ with $\sigma_{a+b}$ we 
may now assume the dimension of $V$ to be $a+2p$.

Let $\Ip(p)$ be the $I(p)$ intersected with $\Hom_k(\wedge^p V^\prime,
\wedge^p V^\prime)$ and let 
\[ A^\prime = k[e_1, \ldots, e_{v-1}]/(e_1^2,\ldots,
e_{v-1}^2) \] 
and so $A_p = A^\prime_p \oplus e_v \cdot A^\prime_{p-1}$.
We are now left with showing that the map
\[ \Ip(p) \oplus \Ip(p-1) \rightarrow \Ip(a+p) \oplus \Ip(a+p-1) \]
which is the same as 
\[ A^\prime_p \oplus e_v A^\prime_{p-1} \mto{\begin{pmatrix}\alpha & 0 \\
\gamma & \beta \end{pmatrix}}
A^\prime_{a+p} \oplus e_v A^\prime_{a+p-1} \]
is injective, where $\alpha = \sigma^\prime_a$, $\beta = \sigma^\prime_a$
and $\gamma = e_v \sigma^\prime_{a-1}$.

Suppose $f \oplus e_v g$ is in the kernel. This means $\beta(e_v g)
= -\gamma(f)$ and $\alpha(f) = 0$.
By induction $\gamma$ is an isomorphism, so this says 
\[ \alpha \circ \gamma^{-1} \circ \beta (e_v g) = 0.\]
But 
\[ \alpha = a! (\sigma^\prime_1)^a, \quad \beta =  a! (\sigma^\prime_1)^a,
\quad \gamma = e_v (a-1)! (\sigma^\prime_1)^{a-1}.\]
Since $\beta$ is injective by induction we get 
$\gamma^{-1} \circ \beta = a \sigma^\prime_1/e_v$ and
$\alpha \circ \gamma^{-1} \circ \beta $ becomes 
$a a! (\sigma^\prime_1)^{a+1}/e_v$ which is an isomorphism by induction on 
$p$. This demonstrates the proposition.
\end{proof}

\subsection{Conditions for $\en \te \alpha$ to map to $\wedge^N V$.}

Let $E$ be the exterior algebra on the vector space $V$ of dimension $v$.
Let $y_0, y_1, \cdots, y_j$ be linearly independent elements in $V$
and $y_J = y_1 \wedge y_2 \wedge \cdots \wedge y_j$ (not including $y_0$).

\begin{lemma} \label{ReLeM} Let $M$ be a homogeneous element in $E$
such that $M \wedge y_J$ has degree $\leq v-2$. If for every $y$ in $V$,
$M \wedge y_J \wedge y$ has $y_0$ as a factor (i.e. can be written
$M^\prime \wedge y_0$ where $M^\prime$ may depend on $y$), then
$M \wedge y_J = N \wedge y_0 \wedge y_J$ for a homogeneous element
$N$.
\end{lemma}

\begin{proof} Let $\{ y_i \} \cup \{ z_k\}$ be a basis for $V$.
Write $M \wedge y_J$ as 
$N \wedge y_0 \wedge y_J + N^{\prime} \wedge y_J$,
where $N^{\prime}$ is an expression in the $z_k$'s.
We then see that $N^{\prime}$ must be zero.
\end{proof}

Now we study maps $\alpha : \wedge^N V \pil V^{\te N-s}$.
We shall also consider the maps from the dual perspective. Let $W = V^*$.
Dualizing $\alpha$ we get linear maps
\[ W^{\te N-s} \mto{\beta} \wedge^N W. \]

\begin{proposition} \label{RePrAi} 
Let $v \geq N+2$. Then the image of 
\[ \wedge^{N+1} V \mto{\en \te \alpha} V^{\te N+1-s} \]
is contained in $\wedge^{N+1-s} V$ if and only if $\alpha$ is $\ent{N-s} \te \Phi$
for some symplectic form $\Phi : \wedge^s V \pil k$.
\end{proposition}

\begin{proof}
The image of $\en \te \alpha$ is in $\wedge^{N+1-s} V$ iff the expression 
\[ x_1 \wedge \beta(x_2, \ldots, x_{N+1-s})\]
is alternating under permutations of $x_1, \ldots, x_{N+1-s}$.
Let $x_2, \ldots, x_{N+1-s}$ be linearly independent.
We can then apply Lemma \ref{ReLeM} as follows. First take $y = x_1$,
$y_0 = x_2$ and $j = 0$. It follows that $\beta(x_2, \ldots, x_{N+1-s})$
contains $x_2$ as a factor. Next let $y=x_1$, $y_0 = x_3$, and $y_1 = x_2$.
It follows that $x_2 \wedge x_3$ is a factor. Thus we may proceed and get
\[ \beta(x_2, \ldots, x_{N+1-s}) 
= x_2 \wedge \ldots \wedge x_{N+1-s} \wedge a \]
for some $a$ (which a priori may depend on the $x$'s).
We also see that
\begin{equation} \beta(x_{\sigma(2)}, \ldots, x_{\sigma(N+1-s)}) =
 x_{\sigma(2)}\wedge \cdots \wedge x_{\sigma(N+1-s)} \wedge a \label{ReLiSi}
\end{equation}
for an arbitrary permutation $\sigma$ of $\{2, \ldots, N+1-s\}$.

If on the other hand $x_2, \ldots, x_{N+1-s}$ span a space of 
dimension $N-s-1$, 
adjoin an $x_{1}$ linearly
independent of these, and suppose if we omit $x_{N+1-s}$
the others span a space of dimension $N-s$. Then switching the first two 
elements
\begin{eqnarray*}
x_1 \wedge \beta(x_2, \ldots, x_{N+1-s})  & = 
& -x_2 \wedge \beta(x_1, \ldots, x_{N+1-s}) \\
& = & -x_2 \wedge x_1 \wedge x_3 \wedge \cdots \wedge x_{N+1-s} \wedge 
a^\prime \\ & = & 0.
\end{eqnarray*}

Since this is true for all $x_{1}$, we must have 
$\beta(x_2, \ldots, x_{N+1-s})$ equal to zero. By descending induction on the
dimension of the space spanned by the arguments, we see that
$\beta(x_2, \ldots, x_{N+1-s})$ is zero when the elements are dependent.

\medskip

Now let $x_1, \ldots, x_v$ be a basis for $W$. Let $I$ and $J$ be disjoint 
subsets of $\{1, \ldots, v\}$ of cardinalities $N-s$ and $s$. Let 
$I$ be $\{i_1, \ldots, i_{N-s} \}$ and $I^\prime$ be 
$\{i_2, \ldots, i_{N-s+1} \}$
where also $I^\prime$ and $J$ are disjoint. We now let
$\beta(x_{i_1}, \ldots ,x_{i_{N-s}})$, which we write $\beta(x_I)$, contain
the term $c_J \cdot x_I \wedge x_J$. Since
\[ x_{i_1} \wedge \beta(x_{i_2}, \ldots, x_{i_{N+1-s}}) 
= (-1)^{N-s} \cdot x_{i_{N+1-s}} \wedge \beta(x_{i_1}, \ldots, x_{i_{N-s}})  \]
we get that
$\beta(x_{I^\prime})$ will contain the term
$c_J \cdot x_{I^\prime} \wedge x_J$.
By using (\ref{ReLiSi}) together with this observation, we easily obtain
\[ \beta(x_I) = c_J \cdot x_I \wedge x_J \]
for arbitrary $I$ which is disjoint from $J$.
Letting $F = \sum_J c_J \cdot x_J$, we have 
$\beta(x_I) = x_I \wedge F$ for all $I$ (where now repetitions in $I$ is 
allowed). Letting $\Phi$ be the dual of $F$ we get $\alpha$ equal to
$\ent{N-s} \te \Phi$.
\end{proof}

\subsection{Conditions for $[\en, \alpha]$ to be zero.}
Now note that when $R$ is $\wedge^N V$ then $(V \te R) \cap (R \te V)$ is
$\wedge^{N+1}V$ by \cite[Lemma 3.12]{Ber01}. 
Also note that if $\alpha_1, \alpha_2, \ldots, \alpha_{i}$ have been found,
we can find a particular solution $\alpha_{i+1}$ to the equation
\[ \alpha_i \circ [\en, \alpha_1] = [\en, \alpha_{i+1}] \]
of (\ref{NKLikJac2}),
and then add solutions of the homogeneous equation 
\[ [\en, \alpha_{i+1}] = 0. \] 
So we want to find conditions on linear maps
\[ \wedge^N V \mto{\alpha} V^{\te N-s}. \]
such that the induced map
\[ \wedge^{N+1} V \mto{ [\en, \alpha]} V^{\te N-s+1}\]
is zero.
Considering this from the dual perspective, we have the dual of 
$\alpha$
\[ W^{\te N-s} \mto{\beta} \wedge^N W, \]
and the dual of $[\en, \alpha]$ is given by the map sending 
$x_1 \cdots x_{N-s+1}$ to 
\[  x_1 \wedge \beta(x_2, \ldots, x_{N-s+1}) -
\beta(x_1, \ldots, x_{N-s}) \wedge x_{N-s+1}.\]

\begin{proposition}\label{RePrKe}
Consider the following maps from $\Hom(\wedge^N V, V^{\te N-s})$ to
$\Hom(\wedge^{N+1} V, V^{\te N+1-s})$

a. \[ \alpha \mapsto \en \te \alpha  + (-1)^s \cdot \alpha \te \en. \]
If either $v \geq N+2$, or $v=N+1$ with $N-s$ even, this map is injective.

b. \[ \alpha \mapsto \en \te \alpha + (-1)^{s-1} \cdot \alpha \te \en. \]
If $v \geq N+2$ then $\alpha$ is in the kernel if and only if
$\alpha$ is $\ent{N-s} \te \Phi$ for a form $\Phi$ in 
$\Hom(\wedge^s V, k)$.
\end{proposition}

\begin{proof}
Let $\delta$ be $0$ or $1$.
The expression $\en \te \alpha  + (-1)^{s-\delta} \cdot \alpha \te \en$ 
is zero iff
\begin{eqnarray} \label{ReLaBe}
& & x_1 \wedge \beta(x_2, \ldots, x_{N+1-s })  \\
\notag &=&
(-1)^{s+1-\delta} \cdot \beta(x_1, \ldots, x_{N-s}) \wedge x_{N+1-s} \\ \notag
& = & (-1)^{N+s+1-\delta} \cdot x_{N+1-s} \wedge \beta(x_1, \ldots, x_{N+1-s}) 
\end{eqnarray} 
Shifting $x_1, \ldots, x_{N-s}, x_{N+1-s}$ cyclically $N+1-s$ times gives
\[ x_1 \wedge \beta(x_2, \ldots, x_{N+1-s}) = (-1)^{(N+1-s)(N+1-s-\delta)}
\cdot x_1 \wedge \beta(x_2, \ldots, x_{N+1-s}). \]
Hence if $\delta = 0$ and $N-s$ is even we get that $\beta$ must be zero.

\medskip Now assume $v \geq N+2$. 
The same argument as in Proposition \ref{RePrAi}
gives (assuming $\alpha$ is in the kernel)
\[ \beta(x_1, \ldots, x_{N-s}) = x_1 \wedge \cdots \wedge x_{N-s} \wedge a \]
for some $a$ (which a priori may depend on the $x$'s). 
In particular 
if $x_i$ and $x_j$ are equal for some $i < j$ then 
$\beta(x_1, \ldots, x_{N-s})$ is zero. But then
$\beta$ factors through $\wedge^{N-s} V$ so
\[ \beta(x_{\sigma(1)}, \ldots, x_{\sigma(N-s)}) = 
x_{\sigma(1)} \wedge \cdots \wedge x_{\sigma(N-s)} \wedge a\]
for any permutation $\sigma$.

Again exactly the same argument as in Proposition \ref{RePrAi} gives
that there exists a form $F$ of degree $s$ such that $\beta(x_I)$ is 
$x_I \wedge F$.
Putting this into (\ref{ReLaBe}) we see that it holds if $\delta$ is
one. If $\delta$ is zero one gets
\[ x_1 \wedge \cdots \wedge x_{N+1-s} \wedge F =
- x_1 \wedge \cdots \wedge x_{N+1-s} \wedge F \]
for all $x_1, \ldots, x_{N+1-s}$, and so in this case $F$ is zero.
\end{proof}

\subsection{Conditions for $[\en, \alpha]$ to map to $\wedge^N V$.}

Now we consider condition (\ref{NKLikJac1}) of Theorem \ref{NKThmPBW}. 
This will separate into two cases.
For naturality of proof we 
formulate slightly more general results than what is needed later on. 
For 
\[ \wedge^N V \mto{\alpha}  V^{\te N-s}\] 
we have the bracket
\[ \{\en,\alpha\} = \en \te \alpha  + (-1)^s \cdot \alpha \te \en \]
which is a linear map from
$\wedge^{N+1} V$ to $V^{\te N+1-s}$. 

\begin{proposition} \label{RePrLi} Suppose $N-s$ is even and $v \geq N+1$.
In the map 
\[ \Hom(\wedge^N V , V^{\te N-s} ) \pil \Hom(\wedge^{N+1} V, V^{\te N+1-s})
 \]
given by sending $\alpha$ to $\{\en, \alpha\}$, the latter
is contained in
$ \Hom(\wedge^{N+1} V, \wedge^{N+1-s} V)$ if and only if
$\alpha$ is $\ent{N-s} \te \Phi$ for a form $\Phi$ in $\Hom(\wedge^s V,k)$.
\end{proposition}

\begin{proposition} \label{RePrSt} 
Suppose $N-s$ is odd, say $N$ is $2p+s+1$, and $v \geq N+2$.
In the map 
\[ \Hom(\wedge^N V , V^{\te N-s} ) \pil 
\Hom(\wedge^{N+1} V, V^{\te N+1-s}) \]
given by sending $\alpha$ to $\{\en,\alpha\}$, the latter is contained in
\begin{equation} 
\Hom(\wedge^{N+1} V, \wedge^{N+1-s} V) \label{AlfaLikIma1}
\end{equation} 
if and only if $\alpha$ is  $\underline{\en^{2p} L}$ for a linear map $L$ in  
$\Hom(\wedge^{s+1} V, V)$.
\end{proposition}

\begin{proof}
We shall only prove the latter proposition since the former is similar
but easier. It is easy to see that if $\alpha$ is $\underline{\en^{2p}L}$,
then $\{\en, \alpha\}$ is in (\ref{AlfaLikIma1}). We therefore
assume that $\{\en, \alpha \}$ is in (\ref{AlfaLikIma1}) and will show that
it is of the form $\underline{\en^{2p}L}$.

{\it Step 1.} Suppose first $v = N+2$. Then (\ref{AlfaLikIma1})
is naturally isomorphic to $\Hom(\wedge^{s+1} V, V)$. By Proposition 
\ref{RePrTa} an isomorphism is given by
\[ \Hom(\wedge^{s+1} V, V) \mto{T_{N-s}} 
\Hom(\wedge^{N+1} V, \wedge^{N+1-s} V). \]
We consider $\wedge^{s+1} V$ as a quotient space of $V^{\te s+1}$ and
therefore this as the domain of $L$. Also consider $\wedge^{N+1} V$ and
$\wedge^{N+1-s} V$ as subspaces of $V^{\te N+1}$ and 
$V^{\te N+1-s}$ respectively. The map $T_{N-s}$ is then given by
\[L \mapsto 1/(s+1)!\cdot 
\{\en, \underline{\en^{2p}L} \}. \]
Thus if $\{\en,\alpha\}$ is in (\ref{AlfaLikIma1}) 
we get a form $L$ in $\Hom(\wedge^{s+1} V, V)$ such that

\[ \{\en, \alpha - 1/(s+1)! \cdot \underline{\en^{2p}L} \} \]
becomes zero. Proposition \ref{RePrKe}.a then gives
$\alpha$ equal to $1/(s+1)! \cdot \underline{\en^{2p}L}$.

\medskip 
{\it Step 2.} Suppose now that $v \geq N+3$ and write 
$V = V^\prime \oplus \langle e \rangle$.
Our approach to show that $\alpha$ is of the form $\underline{\en^{2p}L}$ 
will be to use the decomposition of $V$
to write $\Hom(\wedge^{s+1} V, V)$ as 
\begin{eqnarray} \label{ReLaDwe}
 \Hom(\wedge^{s+1} V^\prime, V^\prime) 
&\oplus & \Hom(\wedge^s V^\prime \wedge e, V^\prime)\\ \notag
\oplus \,\,\,\,
\Hom(\wedge^{s+1} V^\prime , \langle e \rangle) 
&\oplus & \Hom(\wedge^s V^\prime \wedge e, 
\langle e \rangle )
\end{eqnarray}
and find the components of $L$ in each of these summands by induction.

The map $\alpha$ now lives in 
\begin{eqnarray} \label{ReLaDc} &&
\Hom(\wedge^N V^\prime , V^{\prime \te N-s })\\ \notag
& \oplus &  \Hom(\wedge^N V^\prime, 
\underset{a+b = N-s-1}{\bigoplus}
V^{\prime \te a} \te e \te V^{\prime \te b}) \\ \notag 
& \oplus & \Hom(\wedge^{N-1} V^\prime \wedge e , V^{\prime \te N-s}) \\ \notag
& \oplus & \Hom(\wedge^{N-1} V^\prime \wedge e, 
\underset{a+b = N-s-1}{\bigoplus}
V^{\prime \te a} \te  e \te V^{\prime \te b})
\end{eqnarray}
together with $\Hom$-terms where the codomain involves two or more $e$'s,
and is sent to $\{\en, \alpha \}$ which lives in 
\begin{eqnarray}\label{ReLaDc2} &&
\Hom(\wedge^{N+1} V^\prime , V^{\prime \te N+1-s }) \\ \notag &\oplus &  
\Hom(\wedge^{N+1} V^\prime, \underset{a+b = N-s}{\bigoplus}
V^{\prime \te a} \te e \te V^{\prime \te b}) \\ \notag &\oplus &
\Hom(\wedge^{N} V^\prime \wedge  e , V^{\prime \te N+1-s}) \\ \notag &\oplus &
\Hom(\wedge^{N} V^\prime \wedge e, \underset{a+b = N-s}{\bigoplus}
V^{\prime \te a} \te e \te V^{\prime \te b})
\end{eqnarray}
together with $\Hom$-terms whose codomain also involves two or more $e$'s.

We are then assuming $\{\en,\alpha\}$ is in the subspace
\[ \Hom(\wedge^{N+1} (V^\prime \oplus \langle e \rangle), 
\wedge^{N+1-s}(V^\prime \oplus 
\langle e \rangle )) \]
of (\ref{ReLaDc2}).

\medskip 
{\it Step 3.} 
The first summand of (\ref{ReLaDc}) is the only one mapping nonzero
to the first summand in (\ref{ReLaDc2}).
By induction the component of $\alpha$ in the first summand of 
(\ref{ReLaDc}) is given by $\underline{\en^{2p}L^{\prime}}$ for some $L^\prime$
in $\Hom(\wedge^{s+1} V^\prime, V^\prime)$. 
By extending this by zero in (\ref{ReLaDwe}) we can consider it as a map
in $\Hom(\wedge^{s+1} V, V)$ and subtract $\underline{\en^{2p}L^{\prime}}$
from $\alpha$. This new $\alpha$ will map $\wedge^N V^\prime$ to zero in
$V^{\prime \te N-s}$.

\medskip
{\it Step 4.} 
Consider the mapping to the third summand in (\ref{ReLaDc2}). The 
image of the element $\alpha$  comes from the third summand in (\ref{ReLaDc}).
The map between these two identifies as
\[ \Hom(\wedge^{N-1} V^\prime, V^{\prime \te N-s}) \pil 
\Hom(\wedge^N V^\prime, V^{\prime \te N+1-s}) \] 
given by
\[ \gamma \mapsto \en \te \gamma + (-1)^{s-1} \cdot \gamma \te \en. \]
By induction this is in $\Hom(\wedge^N V^\prime, \wedge^{N+1-s} V^\prime)$
iff $\gamma$ is $\underline{\en^{2p}L^{\prime\prime}}$ for some 
$L^{\prime\prime}$ in $\Hom(\wedge^s V^\prime, V^\prime)$. 
Thus letting $\gamma$ be the restriction of $\alpha$
(to the third summand in (\ref{ReLaDc}))
and extending $L^{\prime\prime}$ by zero in (\ref{ReLaDwe})
we get a map in $\Hom(\wedge^{s+1} V, V)$. We subtract 
$\underline{\en^{2p}L^{\prime\prime}}$
from $\alpha$. Now the components of $\alpha$ in the first and third summand
in (\ref{ReLaDc}) become zero.

\medskip 
{\it Step 5.} We now consider the mapping to the second summand in 
(\ref{ReLaDc2}). The image of the element 
$\alpha$ in this summand
comes from the second summand in (\ref{ReLaDc}).
Let $\alpha$ have component $\alpha^i$ in $\Hom(\wedge^N V^\prime, 
V^{\prime \te i} \te e \te V^{\prime \te j})$ where $i$ ranges
from $0$ to $N-s-1$.
The component of $\{ \en, \alpha\}$ in 
$\Hom(\wedge^{N+1} V^\prime, V^{\prime \te a} \te e \te V^{\prime \te b})$ 
is then
\begin{equation} \label{ReLaAla}
\en \te \alpha^{a-1} + (-1)^s \alpha^a \te \en
\end{equation}

Identifying these spaces with 
$\Hom(\wedge^{N+1} V^\prime, V^{\prime \te N-s})$
we wish that all the (\ref{ReLaAla})
be in the subspace $\Hom(\wedge^{N+1} V^\prime, \wedge^{N-s} V^\prime)$ and 
that they
are alternately equal in the sense that for $0 \leq a \leq N-s$
\[ \en \te \alpha^{a}  + (-1)^s \cdot \alpha^{a+1} \te \en =
-(\en \te \alpha^{a-1} 
+ (-1)^s \cdot \alpha^a \te \en) \]
In particular $\en \te \alpha^{N-s-1}$ is in this subspace and so by 
Proposition \ref{RePrAi} we have $\alpha^{N-s-1}$ equal to 
$\ent{N-s-1} \te \Phi^{\prime}$ for
a form $\Phi^\prime$ in $\Hom(\wedge^{s+1} V^\prime, k)$. But then we see 
by Proposition \ref{RePrKe} that
$\alpha^{N-s-j}$ is $\ent{N-s-1} \te \Phi^{\prime}$ for $j$ odd and 
zero when $j$ is even.

Now extend $\Phi^\prime$ by zero in (\ref{ReLaDwe}) to a map in  
$\Hom(\wedge^{s+1} V, V)$.
Note that $\underline{\en^{2p} \Phi^\prime}$ is equal to 
$p+1$ times $\ent{2p} \te \Phi^{\prime}$. 
Replace $\Phi^\prime$ by a multiple and
subtract $\underline{\en^{2p} \Phi^\prime}$
from $\alpha$, giving us that its three first summands
in (\ref{ReLaDc}) are zero.

\medskip
{\it Step 6.} We now consider the mapping to the fourth summand in 
(\ref{ReLaDc2}). Note that the image of the element $\alpha$ now
exclusively comes from the fourth summand in (\ref{ReLaDc}). 
Again let $\alpha$ have component $\alpha^i$ in 
\[ \Hom(\wedge^{N-1} V^\prime
\wedge e, V^{\prime \te i} \te e \te V^{\prime \te j}) \iso 
\Hom(\wedge^{N-1} V^\prime, V^{\prime \te i} \te e \te V^{\prime \te j}). \] 
We see that the
component of $\{ \en, \alpha\}$ in 
\begin{equation} \label{ReLaHi2} \Hom(\wedge^{N} V^\prime \wedge e, 
V^{\prime \te a} \te e \te V^{\prime \te b}) \iso
\Hom(\wedge^{N} V^\prime , V^{\prime \te a} \te e \te V^{\prime \te b}) 
\end{equation}
becomes
\begin{equation} \label{ReLaAla2}
\en \te \alpha^{a-1} - (-1)^s \cdot \alpha^a \te \en . 
\end{equation}
Identifying again the right side of (\ref{ReLaHi2}) with
$\Hom(\wedge^N V^\prime, V^{\prime \te N-s})$ we wish that the
(\ref{ReLaAla2}) shall be in 
$\Hom(\wedge^N V^\prime, \wedge^{N-s} V^\prime)$ and be alternately equal.
In particular $\en \te \alpha^{N-s-1}$ shall be in $\Hom(\wedge^N V^\prime, 
\wedge^{N-s}V^\prime)$.
Proposition \ref{RePrAi} then gives $\alpha^{N-s-1}$ equal to 
$\ent{N-s-1} \te \Phi^{\prime\prime}$ 
for a form $\Phi^{\prime\prime}$ in $\Hom(\wedge^s V^\prime, k)$.
But then we see that $\alpha^{N-s-j}$ is 
$\ent{N-s-1} \te \Phi^{\prime\prime}$ for $j$ 
odd and zero when $j$ is even.

Extend $\Phi^{\prime\prime}$ by zero in (\ref{ReLaDwe}) to a map in 
$\Hom(\wedge^{s+1} V, V)$, replace it by a suitable multiple, and subtract 
$\underline{\en^{2p} \Phi^{\prime\prime}}$ from $\alpha$.
The result is that all components of $\alpha$ in (\ref{ReLaDc}) are zero. 

\medskip
{\it Step 7.}
If some summand of $\alpha$ lives in a $\Hom$-term whose domain is
$\wedge^N V^\prime$ and whose codomain involves two or more $e$'s, then 
our hypothesis on the image $\{\en, \alpha\}$ says that this summand must
map to zero. Using Proposition \ref{RePrAi} we easily deduce that this summand
of $\alpha$ must be zero.

Similarly we can argue for summands of $\alpha$ living in a $\Hom$-term
whose domain is $\wedge^{N-1} V^\prime \wedge e$ and whose codomain
involves two or more $e$'s.
 
The upshot is now that the original $\alpha$ is $\underline{\en^{2p}L}$
where $L$ is $L^\prime + L^{\prime\prime} + \Phi^\prime + 
\Phi^{\prime\prime}$. 
\end{proof}

\section{Auxiliary results}

Now let $A$ be the algebra $T(V) / (xyt-txy)_{x,y,t \in V}$.
Note that $x_I = x_{\sigma I}$ in $A$ for an arbitrary even 
permutation  $\sigma$.
The following will be useful in the proof of Theorem \ref{Se-T2}.

\begin{lemma} \label{TA-T11}
For $p \geq 2$ there is an exact sequence
\[ 0 \pil \wedge^p V \mto{i} A_p \mto{j} S_p(V) \pil 0 \]
where $j$ is the natural quotient map and $i$ is given by
\begin{equation} 
x_{i_1} \wedge \ldots \wedge x_{i_p} \mapsto x_I - x_{\sigma I} \label{TA-2}
\end{equation}
for an arbitrary odd permutation $\sigma$.
\end{lemma}

\begin{proof}
Fix a basis $x_1, \ldots, x_v$ for $V$. Let $x_I$
be a monomial in $A_p$. 
As noted $x_I = x_{\sigma I}$ for any even permutation $\sigma$.
If two indices $i_a$ and $i_b$ in $I$ 
are equal, $\sigma  
I = \sigma (a,b) I$ and so $x_I = x_{\sigma I}$ holds for all permutations.

In any multihomogeneous component of $A$ of degree $(b_1, \ldots, b_p)$ ,
where some $b_i \geq 2$,  $j$ therefore becomes an isomorphism. 
Clearly the elements on the right in (\ref{TA-2}) are non-zero
and linearly independent elements in the kernel when $I$ varies
over all increasing sequences of indices, 
and they also generate the kernel.
\end{proof}

This lemma will be applied as follows. Suppose $N$ is even and 
$\rho : \wedge^N V \pil V^{\te N-s}$ is a map where $s$ is even.
Then $[\en,\rho]$ applied to $x_1 \wedge \cdots \wedge x_{N+1}$ is 
(up to a nonzero constant)
\[ \sum_i x_i \rho(x_{i+1}, \ldots, x_{N+1},x_1,\ldots, x_{i-1})
- \rho(x_{i+1}, \ldots, x_{N+1},x_1,\ldots, x_{i-1})x_i. \]
The crucial thing to note is that this is zero in $A_{N+1-s}$.
Hence a necessary condition for the  equation of Theorem \ref{NKThmPBW}
\[ \alpha_{s-1} \circ [\en,\alpha_1] = [\en,\alpha_s] \]
to be valid is that
the image of the left side is zero in $A_{N+1-s}$.

\medskip
Before taking on the proof of Theorems \ref{Se-T1} and \ref{Se-T2} 
we need some extra lemmas to make the argument go smoother.

In the following let 
\[ L :  \wedge^2 V \pil V, \quad 
\Phi_{2r} : \wedge^{2r} V \pil k\]
 be linear maps. 
We may consider them as maps from  $V^{\te 2}$ and
$V^{\te 2r}$ respectively, by composition with the natural quotient maps.
By combining these maps and the identity maps in various ways we
get maps defined on various $V^{\te q}$ and we shall again restrict these
maps to the subspace $\wedge^q V$ of $V^{\te q}$. 

\begin{lemma} \label{Di-T1}

1.  $L \circ (\en \te L ) = - L \circ (L \te \en)$ on the subspace $\wedge^3 V$
of $V^{\te 3}$.

2. $\Phi_{2r}(\ent{a} \te L \te \ent{b}) = (-1)^b \Phi_{2r}(\ent{a+b} \te L)$
on the subspace $\wedge^{2r+1} V$ of $V^{\te 2r+1}$.

3. Let $M$ be a tensor monomial in $L$ and $\en$ of degree $2r$ with at least 
two $L$'s. If $\Phi_{2r} \circ M$ is defined on $V^{\te q}$, it is zero on
the subspace $\wedge^q V$.
\end{lemma}

\begin{proof}
1. It is immediate to see that this holds for an element $a\wedge b \wedge c$
in $\wedge^3 V$ which is equal to $abc-acb+\cdots$ in $V^{\te 3}$.

2. When $r=1$ we easily check $\Phi_2(\en \te L) = - \Phi_2(L \te \en)$
on $\wedge^3 V$. The general argument is analogous.

3. Note that $\Phi_2(L \te L) = 0$ on $\wedge^4 V$ 
because the image of $L\te L$
is in in the symmetric tensors $S^2(V)$ in $V^{\te 2}$.
The argument in the general case is a straightforward extension.
\end{proof}

\begin{corollary} \label{Di-T41} If
\[ L \te \Phi_{2r}(\underline{\smpm \one^{2r-1} L}) : \wedge^{2r+3} V
\pil V \] is zero, then the following is also the zero map when $b \geq 1$
\[\underline{\en^{2a}L^b} \te
\Phi_{2r}(\underline{\smpm \en^{2r-1} L})
: \wedge^{2(r+a+b)+1} V \pil V^{\te 2a+b}.\]
\end{corollary}

\begin{proof}
Note first that 
\[ \ent{2} \te \Phi_{2r}(\ent{2r-1} \te L) =  
\Phi_{2r} (\ent{2r-1} \te L) \te \ent{2} \] on $\wedge^{2r+3} V$. 
Let $M$ be a tensor monomial in $L$ and $\ent{2}$.
By Lemma \ref{Di-T1} 2. it will be sufficient to show that 
$M \te \Phi_{2r}(\ent{2r-1} \te L)$ becomes zero.
But if $M = M^\prime \te L \te \ent{2b}$ then
\begin{equation} 
M \te \Phi_{2r}(\ent{2r-1} \te L) = M^\prime \te L \te 
\Phi_{2r}(\ent{2r-1} \te L) \te \ent{2b} \label{Di-1}
\end{equation}
For an element $x_1 \wedge \ldots \wedge x_m$ 
(where $m = 2(b+a+r) + 1$)
which we consider as 
\[ \sum (-1)^{sgn(\sigma)} x_{\sigma(1)}\te \cdots \te x_{\sigma(m)} \]
fix the $\sigma(i)$ 
indices which are in positions corresponding to $M^\prime$ and
to $\ent{2b}$ in the right side of (\ref{Di-1}).
By varying the middle terms, the right side of (\ref{Di-1}) 
applied to this becomes zero.
Varying also the end terms, we get the corollary.
\end{proof}

\begin{lemma} \label{DiPrEk} Suppose $L$ is a Lie bracket.
The maps 
\begin{equation} 
L \te \Phi_{2r}(\underline{\smpm \en^{2r-1}L}) : \wedge^{2r+3} V
\pil V \label{Di-3}
\end{equation}
and 
\begin{equation}
L \circ T_2(\Phi_{2r}) \circ T_{2r+1}(L) : \wedge^{2r+3} V \pil \wedge^{2r+2} V
\pil \wedge^2 V \pil V \label{Di-2} 
\end{equation}
are equal up to multiplication by a non-zero scalar.
\end{lemma}

\begin{proof}
The first map in (\ref{Di-2}) identifies up to nonzero scalar as 
$\underline{\smpm \en^{2r+1} L}$
restricted to the subspace $\wedge^{2r+3} V$ of $V^{\te 2r+3}$
The middle map in (\ref{Di-2}) identifies as $\ent{2} \te \Phi_{2r}$
up to a scalar. That the composition in (\ref{Di-2})
now is equal to the one in (\ref{Di-3}) up to a scalar follows since 
$L (\en \te L)$ and $L (L \te \en)$ are zero by the Jacobi identity.
\end{proof}

\begin{lemma} \label{Di-T21}
Suppose $L$ is a Lie bracket. Then
\begin{equation} \underline{\en^{2a} L^b} \circ 
[\en, \underline{\en^{2(a+b-1)}L}] 
= [ \en,\underline{\en^{2(a-1)}L^{b+1}}].\label{Di-11}
\end{equation}
\end{lemma}

\begin{proof}
We use induction on $a+b$.
The first of the terms in (\ref{Di-11}) can be written 
\begin{equation}
\underline{\en^{2a}L^{b}} = \ent{2} \te
\underline{\en^{2(a-1)}L^b} + 
L \te \underline{\en^{2a}L^{b-1}}. \label{Di-12}
\end{equation}
The second of the terms in (\ref{Di-11}) can be written
\begin{equation}
\ent{2} \te [\en, \underline{\en^{2(a+b-2)}L}] + [\en,L] \te \ent{2(a+b-1)}.
\label{Di-13} 
\end{equation}
The terms in (\ref{Di-13}) composed with the first term in 
(\ref{Di-12}) become by induction
\begin{equation}
\ent{2} \te [\en, \underline{\en^{2(a-2)}L^{b+1}}] + [\en,L] \te
\underline{\en^{2(a-1)}L^b}, \label{Di-14}
\end{equation}
and the terms in (\ref{Di-13}) composed with the second term in 
(\ref{Di-12}) become (using that $L \circ [\en,L]$ is zero since $L$ is
a Lie bracket)
\begin{equation}
L \te [\en, \underline{\en^{2(a-1)}L^b}]. \label{Di-15}
\end{equation}

 
Thus the left side of equation (\ref{Di-11}) is
the sum of (\ref{Di-14}) and (\ref{Di-15}) which is 
\begin{eqnarray*}
& \ent{3} \te
\underline{\en^{2(a-2)}L^{b+1}} -
& \ent{2} \te \underline{\en^{2(a-2)}L^{b+1}} \te \en \\
+& \en \te L \te 
\underline{\en^{2(a-1)}L^b} 
-&L \te \underline{\en^{2(a-1)}L^b} \te \en 
\end{eqnarray*}
By using (\ref{Di-12}) this is seen to be equal to the right side
of  (\ref{Di-11}).
\end{proof}

\section{Proof of Theorems \ref{Se-T1} and \ref{Se-T2} }

We are now ready to give the proofs of the main theorems about
PBW-deformations of $T(V)/(\wedge^NV)$.  First we do the case when $N$ is even.

\begin{proof}[Proof of Theorem \ref{Se-T2}]

{\it Step 1.}
By Proposition \ref{RePrSt},
$\alpha_1 = \underline{\en^{2(n-1)}L}$.
When $N=2$ we know $L$ is a Lie bracket. So assume $N\geq 4$.
The composition $\alpha_1 \circ [\en,\alpha_1]$
is then
\[ \underline{\en^{2(n-1)}L } \circ \underline{\smpm \en^{2n-1}L}.\]
We do not yet know that $L$ is a Lie bracket, but the same argument
as in Lemma \ref{Di-T21} shows (or simply check it directly) that
this is
\begin{equation}
[\en, \underline{\en^{2(n-2)}L^2}] + \sum_{a+b=n-1} \ent{2a} \te
L \circ (\en \te L - L \te \en) \te \ent{2b}. \label{BeLaSs}
\end{equation}
Now this is to be equal to $[\en,\alpha_2]$ and so shall be zero in
$A_{N-1}$ of Lemma \ref{TA-T11}. The first term is zero in $A_{N-1}$.
The other term is $0$ in $S(V)_{N-1}$, seen as follows. 
Recall that we consider $[\en,\alpha_1]$ on the subspace $\wedge^{N+1} V
$ of $V^{\te N+1}$, when $a \geq 1$ just keep 
$x_{\sigma(3)} \ldots x_{\sigma(N+1)}$ fixed and switch $x_{\sigma(1)}$ and 
$x_{\sigma(2)}$. Similarly we may argue when $b \geq 1$.
Hence the image in $A_{N-1}$ of the second term lies naturally in 
$\wedge^{N-1} V$
by the sequence of Lemma \ref{TA-T11}.
In fact it is equal to a multiple of 
\[ T_{N-2}(L \circ (\one \te L - L \te \one))\]
which by Lemma \ref{Di-T1} is a multiple of $T_{N-2}(L \circ (\one \te L))$.
Applying Proposition \ref{RePrTa} and noting that $v \geq (N-2) + 3 + 1 $
we must have $L \circ (\en \te L)$ equal to zero if (\ref{BeLaSs}) is zero 
in $A_{N-1}$. This is the Jacobi identity and
so $L$ is a Lie bracket.

\medskip
{\it Step 2.} We shall now proceed by induction on  $2r, 2r+1$.
Suppose $\alpha_{2i}$ has the form (\ref{Se-5}) for $1 \leq i \leq r$ and
that (\ref{SeLaGL}) holds for $i < r$.
By the equivalence of Lemma \ref{DiPrEk} we have
\begin{equation} L \te \Phi_{2i}(\underline{\smpm \en^{2i-1}L}) = 0
\label{SeLaSt2Ek} \end{equation}
for $i < r$.
We shall
show that $\alpha_{2r+1}$ must have the form given.
Since $[\en,\alpha_1]$ is equal to $\underline {\smpm \en^{2n-1}L}$
we get by Lemma \ref{Di-T21} that the composition of
$[\en,\alpha_1]$ and $\alpha_{2r}$ is the sum of terms
\begin{equation} \label{BeLa-1}
\sum_{i=0}^r \, \, [\en,\underline{\en^{2(n-1-2r+i)}L^{2r+1-2i}}]
\te \Phi_{2i}
\end{equation}
and terms
\begin{equation}
\sum_{i = 1}^r \, \, \underline{\en^{2(n-2r+i)}L^{2r-2i}} 
\te \Phi_{2i}(\underline{\smpm \en^{2i-1} L}) \label{Be-2}
\end{equation}
By the induction hypothesis and Corollary \ref{Di-T41} we get that
all terms in (\ref{Be-2}) vanish except when $i=r$ where we get 
\begin{equation}\ent{2n-2r} \te 2r\Phi_{2r}(\ent{2r-1} \te L).\label{BeLaSt21}
\end{equation}

Now let $\alpha_{2r+1}^\prime$ be given by 
(\ref{Se-6}) in the statement of Theorem \ref{Se-T2}.
Since $\en \te \Phi_{2i}$ is equal to $\Phi_{2i} \te \en$ 
on $\wedge^{2i+1}V$ and 
\[\en \te \Phi_{2r}(\ent{2r-1} \te L) = -\Phi_{2r}(\ent{2r-1}\te L) \te \en \] 
on $\wedge^{2r+2} V$
we see that $[\en,\alpha^\prime_{2r+1}]$ is equal to 
the composition of $[\en,\alpha_1]$ and $\alpha_{2r}$,
the sum of (\ref{BeLa-1}) and (\ref{BeLaSt21}).
By Proposition \ref{RePrKe}.a we get 
$\alpha_{2r+1} = \alpha^\prime_{2r+1}$.

\medskip
{\it Step 3.}
Now we assume that $0 \leq r \leq n-1$, that 
$\alpha_{2r+1}$ is given by (\ref{Se-6}), and (\ref{SeLaGL}) holds for
for $i < r$, which is equivalent to (\ref{SeLaSt2Ek}).
Now $[\en,\alpha_1]$ is equal to $\underline {\smpm \en^{2n-1}L}$.
Then the composition of $[\en,\alpha_1]$ and 
$\alpha_{2r+1}$ is equal to 
the sum of terms (using the vanishing (\ref{SeLaSt2Ek}), Corollary \ref{Di-T41}
and Lemma \ref{Di-T1})
\begin{equation}
\sum_{i=0}^r \, \, [\en, \underline{\en^{2(n-2r-2+i)}L^{2r+2-2i}}]
\te \Phi_{2i}
\label{Be-9}
\end{equation}
and terms (using Lemma \ref{Di-T1}.3 and 
the Jacobi identity $L(\en \te L) = 0$)
\begin{equation}
\underline{\en^{2(n-r-1)}L} \te \Phi_{2r}(\underline{\smpm \en^{2r-1}L})
- \underline{\smpm \en^{2n-2r-2}L} \te
r\Phi_{2r}(\ent{2r-1}\te L)
\label{Be-11}
\end{equation}
Again by Lemma \ref{Di-T1}.2 the sum of the terms in (\ref{Be-11}) are
\begin{equation}
(\sum_{a+b = 2n-2r-2} \ent{a} \te L 
\te \ent{b}) \te  
r\Phi_{2r}(\ent{2r-1} \te L). \label{BevLikRPhi}
\end{equation}

We now want to show that 
\begin{equation} L \te \Phi_{2r}(\ent{2r-1} \te L) = 0. \label{BevLikGL}
\end{equation}
which by Lemma \ref{DiPrEk} is the same as (\ref{SeLaGL}).
If the sum of (\ref{Be-9}) and (\ref{Be-11}) 
is to be equal to $[\en,\alpha_{2r+2}]$ then
it has to be zero in $A_{N+1-2r-2}$ of Lemma \ref{TA-T11}.
Since (\ref{Be-9}) vanishes there and (\ref{BevLikRPhi}) becomes zero 
in $S(V)_{N+1-2r-2}$ we must have (\ref{BevLikRPhi})  
zero in $\wedge^{N+1-2r-2} V$. By Proposition \ref{RePrTa}, since $v \geq 
(N -2r-2) + (2r+3) + 1$, we get (\ref{BevLikGL}).
Then we see that (\ref{Se-5}) for 
$\alpha_{2r+2}$ makes $[\en,\alpha_{2r+2}]$
equal to the composition of $[\en, \alpha_1]$ and $\alpha_{2r+1}$ 
and by Proposition \ref{RePrKe}.b every 
solution $\alpha_{2r+2}$ has this form.
\end{proof}

Now we give the much easier proof of the case when $N$ is odd.

\begin{proof}[Proof of Theorem \ref{Se-T1}]
Proposition \ref{RePrLi} gives $\alpha_1$ is $\ent{N-1}\te l$.
Since $l\te \en = -\en \te l$ on $\wedge^2 V$ we get 
$[\en,\alpha_1]$ equal to $\ent{N}\te 2l$.

Now if $\alpha_{2r}$ is given as $\ent{N-2r} \te \Phi_{2r}$ then
$\alpha_{2r} \circ [\en,\alpha_1]$ is 
\begin{equation} \ent{N-2r} \te 2l \te \Phi_{2r} \label{Be-7}
\end{equation}
On the other hand if 
\[ \alpha_{2r+1}^\prime = \ent{N-2r-1} \te l \te \Phi_{2r} \]
then $[\en,\alpha_{2r+1}^\prime]$ is given by (\ref{Be-7}) also.
Proposition \ref{RePrKe} gives $\alpha_{2r+1}$ equal to
$\alpha^\prime_{2r+1}$.

Lastly if $\alpha_{2r+1}$ is given as in (\ref{SeLigAll}), 
then $\alpha_{2r+1} \circ
[\en,\alpha_1] $ is zero since $l\circ(\en\te l)$ vanishes on $\wedge^2 V$.
Hence Proposition \ref{RePrKe}.b gives the form of $\alpha_{2r+2}$.
\end{proof}

\section{PBW-deformations of cubic Artin-Schelter algebras.}

This section is an extended example of Poincar\'{e}-Birkhoff-Witt
deformations of an interesting class of $3$-Koszul algebras; the cubic {\em
Artin-Schelter} algebras.  So we will present equations for these
deformations for the generic forms of these algebras, as presented in
\cite{AS87}.  We give the full details of the argument in one case, a
semi-detailed description in one case,
and as the other cases are similar, we simply present the results.

Throughout this section, the vector space $V$ of generators for our
algebras will be two-dimensional, with basis $x,y$.  The space of
generators $R$ will also be two-dimensional, with basis $f,g$.  The
generators are cubic, so $N=3$.  Finally,
in \cite{AS87} it is proved that in all the situations considered
here, the vector space of syzygies is one-dimensional, and it has a
generator $w$.  We use the conditions from Theorem \ref{NKThmPBW}
throughout.\\

\noindent {\bf Type E.}
From \cite{AS87} we have that this algebra is generated by the two
cubic polynomials 
\[f=y^3+x^3, \quad g=y^2x+\zeta yxy+\zeta^2 xy^2. \] 
Here $\zeta$ is a primitive third root of unity.  In
our notation, $R=(f,g)$, and then $(V\otimes R) \cap (R\otimes
V)$ is 
\[(w) = (y^3x+\zeta y^2xy+\zeta^2 yxy^2+xy^3+x^4).\]  We define general
maps $\alpha_1:R\rightarrow V\otimes V$, $\alpha_2:R\rightarrow V$,  and
$\alpha_3:R\rightarrow k$ with parameters

\[ \begin{array}{l}\alpha_1(f)=a_{11}x\otimes x+a_{12}x\otimes y +
a_{13}y \otimes x + a_{14}y\otimes y \\
\alpha_2(f)=a_{21}x+a_{22}y\\
\alpha_3(f)=a_3, \end{array} \]
and $\alpha_i(g)$ is given similarly with coefficients $b_{ij}$.
Condition (\ref{NKLikJac1}) can be written
\[[\en, \alpha](w)= \beta f + \gamma g,\]
where $\beta$ and $\gamma$ are taken from the ground field.  Now $w$
can be written
\[w=fx+\zeta gy=xf+yg.\]
Thus
\[\beta f +\gamma g= x\alpha_1(f) +  y \alpha_1(g) -\alpha_1(f)x  -\zeta
\alpha_1(g)y. \]

By comparing the coefficients of monomials in $x$ and $y$, we get the
following equations:
\begin{eqnarray}
\label{e1}x^3:  & \beta=a_{11}-a_{11}                       \\
\label{e2}x^2y: & 0= a_{12}-\zeta b_{11}                    \\
\label{e3}xyx:  & 0= a_{13}-a_{12}                          \\
\label{e4}xy^2: & \gamma \zeta^2= a_{14}-\zeta b_{12}       \\
\label{e5}yx^2: & 0=b_{11}-a_{13}                          \\
\label{e6}yxy:  & \gamma \zeta = b_{12}-\zeta b_{13}        \\
\label{e7}y^2x: & \gamma = b_{13}-a_{14}                    \\
\label{e8}y^3:  & \beta=b_{14}-\zeta b_{14}.                 
\end{eqnarray}

Now (\ref{e1}) gives $\beta=0$, and therefore (\ref{e8}) gives $b_{14}=0$.
(\ref{e2}), (\ref{e3}) and (\ref{e5}) give $b_{11}=a_{12}=a_{13}=0$,
whereas (\ref{e4}), (\ref{e6}) and (\ref{e7}) give $b_{12}+b_{13}=0$,
and then the same three equations easily give
$b_{12}=\gamma\zeta/(1+\zeta), b_{13}=-\gamma\zeta/(1+\zeta)$ and $a_{14}=\gamma(1+2\zeta^2)/(\zeta+1)$.\\

Condition (\ref{NKLikJac2}) gives 
$\alpha_1 \circ [\en, \alpha_1]= [\en, \alpha_2]$.  We write
out the coefficients as above, using the fact that $\beta$ has to be
zero: 

\[\gamma \alpha_1(g)=x\alpha_2(f)+y\alpha_2(g)-\alpha_2(f)x-\zeta
\alpha_2(g)y\] 
giving the equations
\begin{eqnarray}
\label{e9}x^2: &0=a_{21}-a_{21}  \\
\label{e10}xy:  &\gamma b_{12}=a_{22}-\zeta b_{21}  \\
\label{e11}yx:  &\gamma b_{13}=b_{21}-a_{22}  \\
\label{e12}y^2: &0=b_{22}-\zeta b_{22}. 
\end{eqnarray}\\

From (\ref{e12}) we get that $b_{22}=0$, from (\ref{e10}) and
(\ref{e11}), using also $b_{12}+b_{13}$, we get $b_{21}=0$.  From
(\ref{e10}) and the computation for $b_{12}$ above, we get
$a_{22}=\gamma^2\zeta/(1+\zeta)$.\\

Condition (\ref{NKLikJac2}) also gives that 
$\alpha_2 \circ [\en, \alpha_1] = [\en, \alpha_3]$.
Going about as above, we find equations

\begin{eqnarray}
\label{e13}x: &0=a_3-a_3  \\
\label{e14}y: &0=b_3-\zeta b_3. 
\end{eqnarray}

This gives $b_3=0$. Finally  $\alpha_3  \circ [\en, \alpha_1]=0$ gives
$0=\beta a_3+\gamma b_3$, which we already knew.\\

There are no restrictions on $a_{11},a_{21}$ and $a_3$,
and all the other coefficients are zero or expressed in terms of
$\gamma$.  To conclude, the deformations we are interested in have two
generators $x$ and $y$, and two relations

\[\begin{array}{l}
y^3+x^3+a_{11}x^2+\frac{\gamma(1+2\zeta^2)}{(1+\zeta)}y^2+a_{21}x+\frac{\gamma^2\zeta}{(1+\zeta)}y+a_3,\\
y^2x+\zeta yxy+\zeta^2 xy^2+\frac{\gamma \zeta}{(1+\zeta)}xy-\frac{\gamma\zeta}{(1+\zeta)}yx\end{array}\]

For the other cases, we will give the coefficients explicitly; then
the expressions for the deformations are given as above.\\

\noindent {\bf Type H.}
With notation as in the previous example, we get from \cite{AS87} that
\[f=-\zeta^3 y^3+\zeta^2 yx^2+\zeta xyx+x^2y, \quad g=y^2x-\zeta yxy
+\zeta^2 xy^2+\zeta^3 x^3. \]  
Here $\zeta$ is a primitive eighth root
of unity.  Then 
\[w=y^3x-\zeta y^2xy+\zeta^2 yxy^2-\zeta^3 xy^3
+x^3y+\zeta x^2yx +\zeta^2 xyx^2 + \zeta^3 yx^3.\]

Define the coefficients of the $\alpha_i$ as before, and also $\beta$
and $\gamma$.  Note that 
\[w=xf+yg=\zeta f x - \zeta gy.\]
From condition (\ref{NKLikJac1}) we get
\[x\alpha_1(f)+y\alpha_1(g)-\zeta\alpha_1(f)x +\zeta\alpha_1(g)=\beta f+\gamma g.\]
Similarly, from (\ref{NKLikJac2}) we get 
\[\begin{array}{l}
\alpha_1(\beta f+\gamma g)=x\alpha_2(f)+y\alpha_2(g)-\zeta\alpha_2(f)x +\zeta\alpha_2(g)y,\\
\alpha_2(\beta f+\gamma g)=x\alpha_3(f)+y\alpha_3(g)-\zeta\alpha_3(f)x +\zeta\alpha_3(g)y,\\
\alpha_3(\beta f+\gamma g)=0.\end{array}\]

As before, by rewriting these equations monomial by monomial, we end
up with expressing the $a$'s and $b$'s in terms of $\beta$ and
$\gamma$.  These expressions are linear for the $a_{1i}$, quadratic
for the $a_{2i}$, cubic for $a_3$ and similarly for the $b$'s. In
addition, there is one quartic relation between $\gamma$ and $\beta$,
which turns out to be trivially fulfilled.

We get the following expressions for the coefficients, expressed in
terms of $\gamma $ and $\beta$ (calculated in Maple, valid in
characteristic different from $2$):
\[\begin{array}{l}
a_{11} = 1/2\gamma(\zeta^3-\zeta^2-\zeta-1)\\
a_{12} = -1/2\beta(3\zeta^3+3\zeta^2+3\zeta+1)\\
a_{13} = -1/2\beta(3\zeta^3+3\zeta^2-\zeta-3)\\
a_{14} = 3/2\gamma(\zeta^3+\zeta^2-\zeta+1)\\
a_{21} =\gamma\beta\zeta(\zeta^2+2\zeta+1)\\
a_{22} =
3/2\beta^2\zeta^3-3/2\zeta\beta^2-2\beta^2+3/2\gamma^2\zeta^3-3\gamma^2\zeta^2+3/2\gamma^2\zeta\\
a_3 = -1/2\gamma(2\beta^2\zeta^3+2\gamma^2\zeta^3+\beta^2\zeta^2-\gamma^2\zeta^2-\zeta\beta^2-\gamma^2\zeta-2\beta^2+2\gamma^2)\\
b_{11} = -3/2\beta(\zeta^3-\zeta^2-\zeta-1)\\
b_{12} = 1/2\gamma(3\zeta^3-3\zeta^2-\zeta+3)\\
b_{13} = 1/2\gamma(3\zeta^3-3\zeta^2+3\zeta-1)\\
b_{14} = -1/2\beta(\zeta^3+\zeta^2-\zeta+1)\\
b_{21} =
-3/2\beta^2\zeta^3-3\beta^2\zeta^2+3/2\gamma^2\zeta-3/2\gamma^2\zeta^3-2\gamma^2-3/2\zeta\beta^2\\
b_{22} = -\gamma\beta\zeta(\zeta^2-2\zeta+1)\\
b_3 = 1/2\beta(2\beta^2\zeta^3+2\gamma^2\zeta^3+\beta^2\zeta^2-\gamma^2\zeta^2-\zeta\beta^2-\gamma^2\zeta-2\beta^2+2\gamma^2)
\end{array}\]

We introduce the other cases by giving the generators $f,g$ of the
ideal, and $w$, a generator for $(V\otimes R) \cap (R\otimes V)$.  These
polynomials we take directly from \cite{AS87}.  Then we give the
conditions on the coefficients of the $a$'s and $b$'s, in terms of
$\gamma$ and $\beta$.  Note that in some cases, there are also
conditions on $\gamma$ and $\beta$.  In a couple of cases, we
have used simplifications like $2\beta=0$ implies $\beta=0$, so there are some
variants in characteristic two that are left out. Many of the
calculations have been carried out using Maple.\\\

\noindent{\bf Type A.}
Here ($a$ and $b$ are general parameters)
\[f=ay^2x+byxy+axy^2+x^3, \quad g=y^3+ayx^2+bxyx+ax^2y\] and thus
\[w=y^4+a(x^2y^2+xy^2x+y^2x^2+yx^2y)+b(xyxy+yxyx)+x^4.\] 
In this case
$\beta=\gamma=0$, and the only conditions on the $a$'s and $b$'s are
\[a_{14}=b_{13}=b_{12},\, a_{13}=a_{12}=b_{11},\, b_{21}=a_{22}.\]
So there are no conditions on the other coefficients.\\

\noindent{\bf Type $\mathbf S_1$.}
In the case $S_1$, we have ($a$ is a general parameter, $\alpha$ is invertible)
\[ f=\alpha xy^2+\alpha^2 y^2x+a\alpha
yxy, \quad g=x^2y+\alpha yx^2+axyx\] and 
\[w=x^2y^2+\alpha
xy^2x+\alpha^2y^2x^2+\alpha yx^2y+axyxy+a\alpha yxyx.\]  We will
consider first the case of generic $\alpha$, meaning different from
$1$, and then the special case afterwards.\\

Generic $\alpha$:
In this case, there is no condition on $\gamma$ or $\beta$, and all
the coefficients can be expressed using these parameters, and in
addition $b_{21}$ (we could also have used $a_{22}$).  The equations are: 

\[\begin{array}{l}
a_{11}=0\\
a_{12}=-\gamma(2+a)\alpha/(\alpha-1)\\
a_{13}=-\gamma(2\alpha+a)\alpha/(\alpha-1)\\
a_{14}= \beta\alpha^2(1+a+\alpha)/(\alpha-1)\\
a_{21}= \gamma^2(1+a+\alpha)\alpha/(\alpha-1)^2\\
a_{22}=\alpha b_{21}\\
a_3=
-\gamma(\beta\gamma\alpha+\beta\gamma\alpha a+\alpha^2\beta\gamma+b_{21}-2b_{21}\alpha+b_{21}\alpha^2)\alpha/(\alpha-1)^3\\
b_{11}= -\gamma(1+a+\alpha)/(\alpha-1)\\
b_{12}= \alpha\beta(2+a)/(\alpha-1)\\
b_{13}= \alpha\beta(2\alpha+a)/(\alpha-1)\\
b_{14}=0\\
b_{22}= \beta^2\alpha^2(1+a+\alpha)/(\alpha-1)^2\\
b_3=
\beta(\alpha b_{21}-2\alpha b_{21}\alpha+\alpha
b_{21}\alpha^2+\alpha^2\beta\gamma+\alpha^2\beta\gamma a+\alpha^3\beta\gamma)/(\alpha-1)^3
\end{array}\]\\

The case $\alpha=1$: We get $\gamma=\beta=0$ (unless $a= -2$).  Then 
\[b_{11}=a_{12}=a_{13},\, b_{12}=b_{13}=a_{14}, \,a_{22}=b_{21}.\]  
There are no other relations, and the other coefficients are free.\\

When $a=-2$ (and still $\alpha=1$), we get a number of conditions as follows: 
\[\begin{array}{l}
a_{12}=b_{11}+\gamma, \, a_{13}=b_{11}-\gamma,\,  b_{12}=a_{14}-\beta,\,  
b_{13}=a_{14}+\beta\\
\beta a_{14}+\gamma b_{14}=0=\beta a_{14}+\gamma b_{14},\,\beta
a_{12}+\gamma b_{12}=0=\beta a_{13}+\gamma b_{13}\\

\beta b_{11}+\gamma a_{14}=0,\, \beta a_{21}+\gamma b_{21}=0=\beta
a_{22}+\gamma b_{22}\\
b_{21}=a_{22}, \, \beta^2 a_{21}=\gamma^2 b_{22},\, \beta a_3+\gamma b_3=0\end{array}\]

\noindent{\bf Type $\mathbf S_2$.}
We have \[f=xy^2+\alpha y^2x, \quad g=-\alpha x^2y+\alpha^2 yx^2\] 
and \[w=x^2y^2+\alpha xy^2x+\alpha^2y^2x^2-\alpha yx^2y.\] As in the previous
case, there are special values of $\alpha$ that require separate
treatment; $\alpha=\pm 1$.\\

First consider the generic case ($\alpha\neq \pm 1$).
Here all the coefficients can be expressed using $\gamma$ and $\beta$,
which are free parameters:\\
\[\begin{array}{llll}
a_{11}=0 &
a_{12}=-{\frac {2\alpha \gamma}{1+\alpha}} &  
a_{13}=-{\frac {2{\alpha}^{2}\gamma}{1+\alpha}}   &
a_{14}={\frac { \left( 1+\alpha \right) \beta}{\alpha-1}}  \\
a_{21}={\frac {{\alpha}^{2}{\gamma}^{2}}{1+\alpha}}  &
a_{22}={\frac {2\alpha\beta\gamma}{1-\alpha}}  &
a_3={\frac {{\alpha}^{2}{\beta}\gamma^{2}}{{\alpha}^{2}-1}}& \\
b_{11}={\frac {{\alpha}^{2}\gamma \left(1-\alpha \right) }{1+\alpha}}  &
b_{12}={\frac {2\beta\alpha}{1-\alpha}}  &
b_{13}={\frac {2\beta{\alpha}^{2}}{\alpha-1}}  &
b_{14}=0\\
b_{21}=-{\frac {2{\alpha}^{2}\beta\gamma}{1+\alpha}}  &
b_{22}={\frac {\alpha{\beta}^{2}}{\alpha-1}}  &
b_3={\frac {{\alpha}^{2}\beta^2\gamma}{1-{\alpha}^{2}}} &
\end{array}\]\\

The case $\alpha=1$: The relations are as follows.
\[\begin{array}{l}\beta=b_{11}=b_{14}=b_{21}=b_{22}=b_3=0,\\
a_{12}=a_{13}=-\gamma,\,
a_{14}=b_{13}=-b_{12},\, a_{22}=\gamma b_{12}.\end{array}\]

The case $\alpha=-1$:
We get 
\[\begin{array}{l}\gamma=a_{11}=a_{14}=a_{21}=a_{22}=a_3=0,\\
b_{12}=b_{13}=-\beta,\, b_{11}=a_{12}=-a_{13},\, b_{21}=\beta
a_{13}.\end{array}\]\\

\noindent{\bf Type $\mathbf S_2'$.}
Now 
\[f=y^2x+xy^2+x^3, \quad g=yx^2-x^2y\]
\[w=x^2y^2+xy^2x+y^2x^2-yx^2y+x^4.\]  We get 
\[\begin{array}{l}\gamma=b_{11}=b_{14}=b_{21}=b_{22}=b_3=0,\\
 a_{12}=a_{13}=-\gamma,\,a_{14}=b_{13}=-b_{12}, a_{22}=\gamma
 b_{12}.\end{array}\]

\end{document}